\documentclass[10pt]{amsart}
\usepackage{amsmath,amssymb,url,bm}
\usepackage[mathscr]{eucal}
\usepackage[hypertex]{hyperref}
\hypersetup{colorlinks=true, urlcolor=blue, citecolor=blue, linkcolor=blue}

\def\mylistparam
 {\setlength{\topsep}{\smallskipamount}%
  \setlength{\itemsep}{\smallskipamount}
  \setlength{\parsep}{0pt}%
  \setlength{\itemindent}{1em}
  \setlength{\leftmargin}{1.75\parindent}%
  \setlength{\labelwidth}{1.75\parindent}%
 }%
\renewenvironment{description}%
 {\begin{list}{}{\mylistparam
        }}
 {\end{list}}

\let\bit=\bibitem
\def\G{Grothendieck} \def\p.{p.\thinspace}  \def\pp.{pp.\thinspace} 
\def\emdash{\unskip\thinspace\penalty10000---\penalty-500\ignorespaces
  \thinspace} 
\def\dotsp{\dots.\,\ignorespaces}
\let\vf=\varphi \let\ve=\varepsilon \let\?=\overline \let\:=\colon
\let\wt=\widetilde \let\wh=\widehat
\let\go=\mathfrak \let\mc=\mathcal \let\I=\mathbf \let\bb=\mathbb
\bmdefine\balpha{\alpha}
\let\into=\hookrightarrow \let\ox=\otimes \let\x=\times
\let\rra=\rightrightarrows \let\xto=\xrightarrow
\let\To=\longrightarrow \let\onto=\twoheadrightarrow
\DeclareMathOperator\Hilb{{\rm Hilb}} \DeclareMathOperator\hilb{\I{Hilb}}
 
\DeclareMathOperator\Pic{{\rm Pic}} \DeclareMathOperator\pic{\I{Pic}}
\DeclareMathOperator\Div{{\rm Div}} \DeclareMathOperator\IDiv{\I{Div}}
\DeclareMathOperator\Hom{{\rm Hom}} \DeclareMathOperator\sHom{{\it Hom}}
\DeclareMathOperator\LS{{\rm LinSys}} \DeclareMathOperator\Ck{{\rm Coker}}
\def\zar{{(\rm Zar)}} \def\et{\text{(\rm \'et)}} \def\fpqc{{(\rm fpqc)}}
\def\smashedlongrightarrow{\setbox0=\hbox{$\longrightarrow$}\ht0=1pt\box0}
\def\risom{\buildrel\sim\over{\smashedlongrightarrow}}

\belowdisplayskip=\abovedisplayskip

\begin{document}

\title{The Picard Scheme}
\author{Steven L. Kleiman}
 \thanks{It is a great pleasure to thank my good friend Luc Illusie for
carefully reading an earlier version of this paper, and making numerous
comments, which led to many significant improvements.  It is a pleasure
as well to thank Michael Artin, Arthur Mattuck, Barry Mazur, David
Mumford, Brian Osserman, Michel Raynaud, and the referee for their
comments, which led to many other important improvements.}
 \address{Department of Mathematics, MIT, Cambridge, MA 02139, USA}
 \email{kleiman@math.mit.edu}
 \begin{abstract}
This article introduces, informally, the substance and the spirit of
\G's theory of the Picard scheme, highlighting its elegant simplicity,
natural generality, and ingenious originality against the larger
historical record.
 \end{abstract}
\maketitle

\section{Introduction}
\begin{flushright}{\it
 A scientific biography should be written\\
 in which we indicate the ``flow'' of mathematics$\,\dots\mspace{-0.25mu}$\\
 discussing a certain aspect of \G's work, indicating possible roots,\\
 then describing the leap \G\ made from those roots to general ideas,\\
 and finally setting forth the impact of those ideas.}\\
 Frans  Oort \cite[\p.2]{Oh}
\end{flushright}
\bigskip

Alexander \G\ sketched his proof of the existence of the Picard scheme
in his February 1962 Bourbaki talk.  Then, in his May 1962 Bourbaki
talk, he sketched his proofs of various general properties of the
scheme.  Shortly afterwards, these two talks were reprinted in
\cite{FGA}, commonly known as {\bf FGA}, along with his commentaries,
which included statements of nine finiteness theorems that refine the
single finiteness theorem in his May talk and answer several related
questions.

However, \G\ had already defined the Picard scheme, via the functor it
represents, on \pp.195-15,\,16 of his February 1960 Bourbaki talk.
Furthermore, on \p.212-01 of his February 1961 Bourbaki talk, he had
announced that the scheme can be constructed by combining results on
quotients sketched in that talk along with results on the Hilbert scheme
to be sketched in his forthcoming May 1961 Bourbaki talk.  Those three
talks plus three earlier talks, which prepare the way, were also
reprinted in \cite{FGA}.

Moreover, \G\ noted in \cite[\p.C-01]{FGA} that, during the fall of 1961,
he had discussed his theory of the Picard scheme in some detail at
Harvard in his term-long seminar, which David Mumford and John Tate
continued in the spring.  In November 2003, Mumford kindly lent me his
own folder of notes from talks given by each of the three, and notes
written by each of them.  Virtually all the content was published long
ago.

 Those notes contain a rudimentary form of the tool now known
as {\it Castelnuovo--Mumford regularity}.  \G\ mentions this tool in his
commentaries \cite[\p.C-10]{FGA}, praising it as the basis for an
``extremely simple'' proof of a bit weaker version of his third
finiteness theorem.  Mumford sharpened the tool in his book
\cite[Lect.~14]{Mcs}, so that it yields the finiteness of the open
subscheme of the Hilbert scheme that parameterizes all closed subschemes
with given Hilbert polynomial.  

\G\ \cite[\p.221-1]{FGA} correctly foresaw that the Hilbert scheme is
``destined to replace'' Chow coordinates.  As he \cite[\p.195-14]{FGA}
put it, they are ``irremediably insufficient,'' because they ``destroy
the nilpotent elements in parameter varieties.''  Nevertheless, he
\cite[\p.221-7]{FGA} had to appeal to the theory of Chow coordinates to
prove the finiteness of the Hilbert scheme.  So after he received a
prepublication edition of \cite{Mcs}, he wrote a letter on 31 August
1964 to Mumford in which he \cite[\p.692]{Msp} praised the theory in
Lecture 14 as ``a significant amelioration'' of his own.

Mumford made use of the finiteness of the Hilbert scheme in his
construction of the Picard scheme over an algebraically closed field in
\cite[Lect.~19]{Mcs}, whereas Grothendieck took care to separate
existence from finiteness, giving an example in \cite[Rem.~3.3,
\p.232-07]{FGA} over a base curve of a Picard scheme with connected
components that are not of finite type.

Mumford's book \cite{Mcs} was based closely on the lovely course he gave
at Harvard in the spring of 1964.  It was by far the most important
course I ever took, due to the knowledge it gave me and the doors it
opened for me.  During the academic year of 1966--67, I was a Postdoc
under \G\ at the IHES (Institut des Hautes \'Etudes Scientifiques).
When he learned from me that I had taken that course and had advanced
some of the finiteness theory in my thesis \cite[Ch.~II]{TNTA}, he asked
me to write up proofs of his nine finiteness theorems for SGA6
\cite[Exps.~XII, XIII]{SGA6}.

\G, perhaps, figured that I had learned how to prove his nine theorems
at Harvard, but in fact I had not even heard of them.  At any rate, he
told me very little about his original proofs, and left me to devise my
own, which I was happy to do.  There is one exception: the first
theorem, which concerns generic relative representability of the Picard
scheme.  Its proof has a very different flavor, as it involves nonflat
descent, Oort d\'evissage, and representability of unramified functors.
\G\ asked Michel Raynaud to lecture on this theorem and to send me his
lecture notes, which I wrote up in \cite[Exp.~XII]{SGA6}.

My experience led me to study Grothendieck's construction of the Picard
scheme, and to teach the whole theory a number of times.  Further, in
collaboration with Allen Altman, Mathieu Gagn\'e, Eduardo Esteves, Tony
Iarrabino, and Hans Kleppe, I extended some of \G's theory to the
compactified Picard scheme.  The underlying variety had been introduced
via Geometric Invariant Theory in 1964 by Alan Mayer and Mumford in
\cite[\S\,4]{MuWH}.  The scheme has been studied and used by many others
ever since then.

Thus my experience is like the experiences of Nicholas Katz and Barry
Mazur, which were described by Allyn Jackson in
\cite[\p.1054]{NAMS51.9}.  Katz said that \G\ assigned him the topic of
Lefschetz pencils, which was new to him, but ``he learned a tremendous
amount from it, and it had a big effect on my future.''  Mazur said that
\G\ asked him this question posed earlier by Gerard Washnitzer: ``Can the
topology of algebraic variety vary with the complex embedding of its
field of definition?''  Mazur, then a differential topologist, added,
``But for me, it was precisely the right kind of motivation to get me to
begin to think about algebra.''

Both Katz and Mazur then confirmed my impression that our experiences were
typical.  Jackson quotes Katz as saying that \G\  got visitors
interested in something, but with ``a kind of amazing insight into
what was a good problem to give to that particular person to think
about.  And he was somehow mathematically incredibly charismatic, so
that it seemed like people felt it was almost a privilege to be asked to
do something that was part of \G's long range vision of the future.''
Similarly, Mazur said that \G\ had an instinct for ``matching people
with open problems.  He would size you up and pose a problem that would
be just the thing to illuminate the world for you.  It's a mode of
perceptiveness that's quite wonderful and rare.''

I spent the summer of 1968 at the IHES.  \G\ invited me to his home in
Massy-Verri\`eres to discuss my drafts for my contributions to SGA\,6.
His comments ranged from providing insight into the theory of bounded
families of sheaves to criticizing my starting sentences with symbols.%
 \footnote{Many years later, Jean-Pierre Serre told me that he had taught \G\ not
to start sentences with symbols.}
  Again, my experience was typical:%
 \footnote{\G\ gave me another project during my Postdoc.  On
April 18 and 25 that year, he talked in his seminar at the IHES on his
Standard Conjectures and Theory of Motives.  He asked me to write up his
talks, gave me copies of his notes on related matters, and invited me to
his home a year later, in the summer of 1968, to discuss my draft.  That
work too led me to learn some good mathematics and to write several
articles, although they are more expository.  Also, it led to my
co-chairing an organizing committee for an AMS summer research
conference in 1991.

However, my experience was the exception that proved the rule: \G\ had
already asked others to write up his talks; they tried, and gave up!
Also, curiously he never told me about his talk on the Standard
Conjectures at a conference in Bombay, India, in January 1968, let alone
offer me his notes.  Moreover, in the conference
proceedings, his writeup cites a
talk of mine at the IHES, which I never gave, crediting me for an
observation; but it is due to Saul Lubkin, and credited to him in
my writeup \cite[\p.361]{KacWc}, which \G\ critiqued in his home that
summer.}
 Jackson \cite[\p.1054]{NAMS51.9} quotes Luc Illusie as saying, that \G\
often worked at home with colleagues and students, making a wide range
of apposite comments on their manuscripts.

One time, \G\ found that I didn't know some result treated in EGA
(\cite{EGAI} and sequels).  So he gently advised me, for my own good, to
read a little EGA every day, in order to familiarize myself with its
content.  After all, he pointed out, he had been writing EGA as a
service to people like me; now it was up to us to take advantage of this
resource.  That experience supports a statement Leila Schneps made in
\cite[\p.16]{SchnepsCh1}: ``The foundational work that \G\ and [Jean]
Dieudonn\'e were undertaking [was] in the service of all mathematicians,
of mathematics itself.  The strong sense of duty and public service was
felt by everyone around \G.''

As \G\ stated on \p.6 of \cite{EGAI}, he planned to develop in EGA the
ideas he sketched in \cite{FGA}.  He did not succeed.  Nevertheless,
those ideas have become a basic part of Algebraic Geometry.  So they
were chosen as the subject of a summer school held 7--18 July 2003 at
the  ITCP (International Center for Theoretical Physics) in Trieste,
Italy.  The first Bourbaki talk reprinted in \cite{FAG} was not covered;
it treats \G's generalization of Serre duality for coherent sheaves, so
is somewhat apart and was already amply developed in the literature.

The lectures were written up, and published in \cite{FAG}.  As stated on
\p.viii, ``this book fills in \G's outline.  Furthermore, it introduces
newer ideas whenever they promote understanding, and it draws
connections to subsequent developments.''  In particular, I wrote about
the Picard scheme, beginning with a 14-page historical introduction,
which served as a first draft for the present article.

Mumford stated the goal of his book \cite{Mcs} on \pp.vii--viii: ``a
complete clarification of\,\dots the so-called%
 \footnote{On \p.8, the theorem is formulated as ``problem (B),'' and
two analytic solutions are outlined.   On
\p.157, a more precise version is formulated as the ``Fundamental Theorem,''
and given its first algebraic proof.   On \p.169, an important special case
is proved, following \G's somewhat different algebraic treatment.  However,
none of those is called the ``Theorem of Completeness.''}
  Completeness of the Characteristic Linear System of a good complete
algebraic system of curves on a surface\dotsp Until about 1960, no
algebraic proof of this purely algebraic theorem was known\dotsp [Then]
a truly amazing development occurred:'' by combining his results on the
Hilbert scheme and the Picard scheme with Cartier's result, ``that group
schemes in characteristic 0 are reduced,'' \G\ \cite[\pp.221-23,24]{FGA}
obtained in February 1961 an enlightening, purely algebraic proof.
``The key\dots is the systematic use of nilpotent elements.''

\G\ had, moreover, reversed history: he proved Completeness via the
Picard scheme.  By contrast, in December 1904 Federigo Enriques and
sometime in 1905 Francesco Severi gave algebraic proofs of Completeness
from scratch.  In the first half of 1905, on the basis of Enriques's
work, Guido Castelnuovo {\bf introduced} the Picard variety in order
to prove the Fundamental Theorem of Irregular Surfaces.  It asserts the
surprising equality of the four basic invariants: the dimension of the
Picard variety, the irregularity, the number of independent Picard
integrals of the first kind, and half the first Betti number.  \G's
theory, without reference to Completeness, also yields the first part of
the Fundamental Theorem, that the dimension of the Picard variety is
equal to the irregularity in characteristic 0.

Both Enriques's and Severi's proofs have serious gaps, as Severi himself
noted in 1921.%
 \footnote{In 1949, Severi \cite[\p.40]{Sgai} lamented the fact that
``this annoying episode was taken as an article of indictment
for the [crime of] lack of rigor in Italian algebraic geometry!''}
  Severi then proved a more restricted version of Completeness, but one
sufficient for Castelnuovo's work.  Severi's proof was based on Henri
Poincar\'e's construction of a key system of curves.  That construction
appeared in 1910, 1911; it is rigorous, but analytic.  After 1921,
finding a fully rigorous, purely algebraic proof of a suitable version
of Completeness became
a major endeavor \emdash undertaken by Enriques, Severi, and
others\emdash until \G\ finally settled the matter.  Section~2 explains
more fully the history and meaning of Completeness and of the
Fundamental Theorem; Section~5 elaborates on \G's proof.

When \G\ worked on his theory of the Picard scheme, the general
algebro-geometric theory of the Picard variety had been under active
development for nearly fifteen years.  More than twenty mathematicians
had worked on various aspects.  \G\ clarified and settled a number of
issues.  Section~3 explains those issues in chronological order.
Sections~4 and 5 give more detail about \G's advances, which involve
many great innovations.

One issue was a topic of conversation between \G\ and Jacob Murre
sometime in the academic year 1960/61.  Murre told Schneps about
it, and
she \cite[\pp.1--2]{SchnepsCh1} quoted him as saying, ``A very important
unsolved question\,\dots\,{}[was] the behavior of the Picard variety if the
original variety\,\dots\,moved in a system and moreover\emdash and
worse\emdash in[to] characteristic $p>0$\dotsp\penalty-100 I asked
\G\ whether he could explain this behavior\dotsp He said he would
certainly [do so]\dotsp{} Then, in 1962, \G\ completely solved the
question\dotsp{}I attended his Bourbaki lectures, and needless to say, I
was very impressed!''

As it happens, much earlier, in his 1958  talk \cite[\p.118]{Gicm58} at
the  ICM (International Congress of Mathematicians,
\G\ said, ``We shall not give here the precise definition of a `relative
Picard schema', but\,\dots\,if this schema exists then it behaves in the
simplest conceivable way with respect to change of base-space.''  In his
February 1960 Bourbaki talk \cite[\pp.195-15,\,16]{FGA}, he added that
``in particular, the Picard schemes of the fibers'' of a system are the
fibers of the relative Picard scheme, once existence is proved.  In his
February 1961 talk \cite[\pp.212-01]{FGA}, as noted above, he announced
his proof of existence.

Thus when \G\ had succeeded in settling a major issue, such as the
Behavior of the Picard Variety in a Family or the Completeness of the
Characteristic System, he noted the advance, but did not tout it.
Cartier \cite[\p.17]{Ch} describes \G's philosophy as follows: ``\G\ was
convinced that if one has a sufficiently unifying vision of mathematics,
if one can sufficiently penetrate the essence of mathematics and the
strategies of its concepts, then particular problems are nothing but a
test; they do not need to be solved for their own sake.''

One beautiful illustration of \G's ``unifying vision'' is provided by
his theory of the Picard functor.  It is the {\it functor of points\/}
of the Picard scheme \emdash that is, the functor whose values are the
sets of maps from a variable source into the scheme.  Often, a functor
of points is said to provide nothing more than another way of expressing
the universal property of a fine moduli scheme.  That statement is true
for the Hilbert scheme, but a half-truth for the Picard scheme.

What is the universal property of the Picard scheme?  The naive answer
falls short!  However, \G\ saw the hidden common thread in descent of
the base field, Galois cohomology, and sheaf theory; he concluded that
any functor of points must be a sheaf for the fpqc \G\ topology.  Thus
the right Picard functor has to be the sheaf associated to the naive
Picard functor, regarded as a presheaf.  More work with the functor
leads to the construction of the Picard scheme.  It is automatically
compatible with base change, because the Picard functor is so.  Sections 4
and 5 explain all that theory.

In short, Section 2 gives a historical introduction to two venerable
theorems: the Theorem of Completeness of the Characteristic System, and
the Fundamental Theorem of Irregular Surfaces.  Section 3 gives a
historical introduction to the inadequate algebro-geometric theory of
the Picard variety.  Please note: these two introductions are {\bf not}
meant to be either serious historical studies or rigorous mathematical
surveys, but simply fascinating informal accounts, providing background
material for comprehending the nature and extent of \G's advances.

Section 4 explains \G's innovative theory of the Picard functor,
culminating in his main construction of the Picard scheme.  Finally,
Section~5 explains how the theory of the Picard scheme enabled \G\ and
others to provide enlightening treatments of the issues discussed in
Sections 2 and 3.  The discussions in Section 4 and 5 are mathematically
rigorous, but just introductory.  Sources for more information are given
at the beginning of each of Sections 2--5.

There are three minor mathematical novelties below: (1) the proof on
\p.\pageref{eqpa} of the equivalence of the 19th century definition of
the arithmetic genus of a surface and the modern definition, (2) the
algebro-geometric treatment on \p.\pageref{S21C} of Severi's 1921
version of Completeness, and (3) the ``nearly formal'' treatment on
\p.\pageref{GAlb} of the Albanese variety, including duality,
intriguingly announced by Grothendieck on \p.232-14 of his February 1962
Bourbaki talk  \cite{FGA}.

  \newpage
 \section{Irregular Surfaces}
 \begin{flushright}{\it
 But to demonstrate the power of modern abstract ideas\\
 to solve older very concrete problems,\\
 I think that this example is unmatched.}\\
 David  Mumford \cite[\p.7]{Mh}
\end{flushright}
\medskip

In the quotation above, {\it this example\/} refers to \G's treatment of
the Theorem of Completeness of the Characteristic System.  In fact, the
example is the centerpiece of Mumford's article \cite{Mh} in this
volume.  Moreover, Mumford notes that Completeness yields the
Fundamental Theorem of Irregular Surfaces.  Thus if we are to appreciate
the full significance of \G's contribution, then we must review the
history of those two main theorems.  We do so in this section.  First we
pursue, intuitively, the spirit of the original work.  Then we treat
that work rigorously, beginning at the end of this section and
continuing in Section 5.

A number of historical reviews are already available, and served as a
basis for the account here.  Notably, in 1906, Castel\-nuovo and
Enriques wrote%
 \footnote{Please also see their encyclopedia article \cite{CEe} and
Castelnuovo's historical note \cite[\pp.339--353]{CEsa}.}
  one \cite{CE1906} at the request of Emile Picard to be an appendix to
Tome II of his book \cite{PSfa} with Georges Simart.  In 1934, Oscar
 Zariski reviewed various
aspects of the development in different places in his celebrated book
\cite{Zas}.  Those reviews are fairly technical.  In 1994, Fabio
Bardelli \cite{Bcort} wrote a more informal review of the developments
through 1934.  In 1974, Dieudonn\'e published a masterful history of
algebraic geometry, which touches on these theorems in particular.  The
book was translated as \cite{Dhdag} by Judith Sally, and supplemented
with an extensive annotated bibliography.

In 2011, Mumford \cite{Mnams} carefully analyzed the mathematics in a
1936 paper by Enriques on Completeness.  Mumford, in his introduction,
stated his conclusion: ``Enriques must be credited with a nearly
complete [algebro-]geometric proof using, as did \G, higher order
infinitesimal deformations.\,\ldots Let's be careful: he certainly had
the correct ideas about infinitesimal geometry, though he had no idea at
all about how to make precise definitions.''%
 \footnote{Mumford elaborated in his resume\'e at the end: ``Although
[Enriques] gave [infinitesimal deformations] names, they remained in
limbo, without substance, because he did not think of what it meant to
have a function on them.  Grothendieck realized that functions on such
objects should be rings with nilpotent elements, and this gave life to
these infinitesimal deformations."}
   Mumford's article is
preceded by a lovely article by Donald Babbitt and Judith Goodstein
\cite{B-G11}, which focuses on the times, lives, and personalities of
Enriques and his colleagues; please also see their related articles
\cite{B-G09} and \cite{G-B12}.  All the articles mentioned above give
many precise references, which are not repeated here.

Around 1865, Alfred Clebsch caused a sea change in algebraic geometry,
turning it away from the concrete study of particular curves and
surfaces, and toward the abstract study of their {\it birational
invariants\/} \emdash the numbers that depend only on their field of
rational (or global meromorphic) functions.

In 1868, Clebsch considered a connected {\it smooth}%
 \footnote{Also called {\it nonsingular,}  $\wt X$ is defined by polynomials
with Jacobian matrix  of maximal rank.}
  complex projective algebraic surface $\wt X$ of large degree $n$.  He
studied it via its general projection in 3-space, which is a surface
\abovedisplayshortskip 0.0pt plus 4.19998pt minus 2.0998pt
 $$X:f(x,y,z)=0\quad\text{and}\quad n:=\deg f$$
with ``ordinary'' singularities, none at infinity, and no point
at infinity on the $z$-axis.

Clebsch found  the algebraic double integrals on $\wt X$ of the
{\it first kind}\emdash that is, those finite on any bounded analytic
domain of integration\emdash to be of the form
\abovedisplayskip=4.19998pt plus 4.19998pt minus 0.5pt
\belowdisplayskip=\abovedisplayskip
	$$\int\!\!\int\frac{h(x,y,z)}{\partial f/\partial z}\,dx\,dy$$
 where $h$ is a polynomial of degree at most $n-4$ vanishing on the
singular locus,
  $$\Gamma:f,\,\partial f/\partial x,\, \partial f/\partial z,\,
        \partial f/\partial z=0,$$
 a curve of double points.  The number of linearly independent such
integrals became known as the {\it geometric genus} and denoted by
$p_g$.

Clebsch asserted without proof that $p_g$ is a birational invariant.  In
1870, his student, Max Noether, gave an algebraic proof.  In 1869, Arthur
Cayley worked out a formula for the number of independent $h$;
essentially, he found an explicit expression for $F(n-4)$ where $F$ is
the Hilbert polynomial of the homogeneous ideal of $\Gamma$.  The value
$F(n-4)$ was later called the {\it arithmetic genus\/} and denoted by $p_a$.%
 \footnote{Cayley \cite{Cdcs} denoted it by $D$, and called it the {\it
deficiency}.  Picard and Simart \cite[\p.88]{PSfa} denoted it by $p_n$,
and called it the {\it numerical genus}.  Those definitions soon fell
into disuse.}

In 1871, Hieronymous Zeuthen used Cayley's formula to prove
algebraically that $p_a$ too is a birational invariant.  Also in 1871,
Cayley observed that, if $X$ is a ruled surface with plane section of
genus $g$, then $p_a=-g\le0$, although $p_g=0$.

The disagreement between $p_g$ and $p_a$ came as a surprise.  In 1875,
Noether explained it: $F(n-4)$ is the number of independent $h$ only if
$n$ is suitably large.  In any case, $p_g\ge p_a$.  Moreover, if $X$ is
smooth or rational, then $p_g=p_a$.  It was thought that, as a rule,
$p_g=p_a$, and when so, $X$ was dubbed {\it regular}.  The failure of
$X$ to be regular is quantified by the difference $p_g-p_a$; so it
became known as the {\it irregularity}.  Zariski \cite[\p.75]{Zas}
denoted it by $q$; following suit, set
 $$q:=p_g-p_a.$$

In 1884, Picard initiated the study of algebraic simple integrals
	$$\int P(x,y,z)\,dx + Q(x,y,z)\,dy$$
 that are {\it closed}, or $\partial P/\partial y = \partial Q/\partial
x$; they became known as {\it Picard integrals}.  He proved that there
are only finitely many independent such integrals of the {\it first
kind}, those finite on any bounded analytic path of integration; use%
 \footnote{Castelnuovo and Enriques \cite[\p.495]{CE1906} used
$q$, whereas Zariski \cite[\p.162]{Zas} used $r_0$.}
  $s$ to denote their number. Picard noted that, if $X$ is smooth, then
$s=0$.

In 1894, Georges Humbert considered an {\it algebraic system,} or {\it
algebraic family,} of curves.  Its members are the zeros on $X$ of a
polynomial
$$\vf(x,y,z;\lambda_0,\dotsc,\lambda_t)$$
 where the $\lambda_i$ satisfy polynomial equations, which define the
{\it parameter variety} $\Lambda$.  The curves can all contain
common subcurves; some of them are included as {\it fixed components\/}
of the system, and the others, omitted.  The system is said to
be {\it linear\/} if there are homogeneous polynomials $\vf_i(x,y,z)$ of
the same degree with $$\vf = \lambda_0\vf_0+\dotsb+\lambda_t\vf_t.$$

Humbert proved a remarkable result: if $s=0$, then every algebraic
system is a subsystem of a linear system.  That result inspired
Castelnuovo to prove in 1896 that, if $q=0$, then again every algebraic
system is a subsystem of a linear system under a certain restriction,
which Enriques removed in 1899.

In 1897, Castelnuovo fixed a linear system of curves on $X$.  Let $r$ be
its {\it dimension,} the number of linearly independent restrictions
$\vf_i|X$ diminished by 1.  Castelnuovo studied its {\it
characteristic\/} linear system, the system cut out by the $\vf_i$ on a
general member curve $D_\eta$, assuming $D_\eta$ is {\it irreducible},
that is, not the union of two smaller curves.  The characteristic system
has dimension has $r-1$.

Castelnuovo formed the {\it complete,} or largest, linear system on
$D_\eta$ containing the characteristic system.  Let $\delta$ be the
amount, termed the {\it deficiency,} by which the dimension of the
characteristic system falls short of the dimension of its complete
linear system.  Castelnuovo proved that
 \begin{equation}\label{eq<=q} \delta\le q,\end{equation}
 with equality if the linear system consists of all {\it hypersurface
sections\/} of high degree, namely if the $\vf_i(x,y,z)$ generate all
homogeneous polynomials of that degree.

In February 1904, Severi extended Castelnuovo's work.  Severi fixed an
algebraic system of curves on $X$, and a general member $D_\eta$.  He
assumed that $D_\eta$ is irreducible and that $D_\lambda\neq D_\mu$ for
all distinct $\lambda,\,\mu\in\Lambda$.  As $\lambda$ approaches $\eta$
along a path in $\Lambda$, the intersections $D_{\lambda}\cap D_\eta$
approach a limit, which depends only the tangent vector at $\eta$ to the
path.  The various limits form a linear system on $D_\eta$,
parameterized by the projectivized tangent space to
$\Lambda$ at $\eta$.  Thus Severi constructed the {\it characteristic
linear system\/} of the algebraic system.  Set $R:=\dim\Lambda$.  Then this
characteristic system is of dimension $R-1$.

In the algebraic system, form the largest linear subsystem containing
$D_\eta$.  Denote its dimension by $r$.  Form its characteristic linear
system.  Let $\delta$ be its deficiency.  Then its complete linear
system has dimension $r-1+\delta$, and it also contains the
characteristic system of the algebraic system.   Thus
Severi proved that
\begin{equation}\label{eqCpltnss}
R\le r + \delta,
\end{equation}
 with equality if and only if the latter characteristic
system is complete.

In December 1904 Enriques and sometime in 1905 Severi each constructed
an algebraic system with $R=r+q$.  Both constructions are short and
delicate.  Both rely on the completeness of the characteristic system of
certain%
 \footnote{Severi \cite[\p.41]{Sgai} noted that both he and Enriques
believed at the time that they had proved every complete algebraic
system with irreducible general member has a complete characteristic
system!}
  algebraic systems.  Both are flawed, as Severi himself pointed out in
1921.  In 1934, Zariski \cite[\pp.99--102]{Zas} reviewed those
constructions, ``in order to analyze the assumption on which they are
based and for which as yet an algebro-geometric proof is not
available.''

In 1910 and 1911 using a new method of ``normal functions,'' Poincar\'e
gave a rigorous analytic construction of an algebraic system with $r=0$
and $R=q$.  His construction was simplified and developed by Severi in
1921 and Solomon Lefschetz in 1921 and 1924.  In 1934, Zariski
\cite[\pp.169--173]{Zas} reviewed that work too.%
 \footnote{Please also see the
reviews of Mumford \cite[\pp.9--10]{Mcs} and Dieudonn\'e
\cite[\p.53]{Dhdag}.}
He \cite[\p.102]{Zas} noted that ``the value of the construction
of such a system is greater than that of mere example; indeed it is an
essential step in the theory.''

Zariski then derived Severi's May 1905 theorem%
 \footnote{In December 1904, Enriques proved the theorem under the more
stringent, but still sufficient, hypothesis that $D_\eta$ is
``regular,'' later renamed ``regular and nonspecial.''  Please see
Fn.\ref{reg} on \p.\pageref{reg}.}
 that, if there is one system with $R=r+q$, then $R=r+q$ holds for every
complete system whose general member $D_\eta$ is {\it arithmetically
effective\/}; namely, a certain lower semi-continuous combination of its
numerical characters is nonnegative, a common condition (please see
\p.\pageref{ae}).  Hence, by \eqref{eq<=q} and \eqref{eqCpltnss}, if
$D_\eta$ is irreducible too, then its characteristic system is complete
and $\delta=q$.  By Bertini's Theorem, usually $D_\eta$ is irreducible.

The {\bf Theorem of Completeness} came%
 \footnote{According to Severi \cite[\p.42]{Sgai}, in 1921 he derived
the theorem essentially in this form from Poincar\'e's construction.}
  to mean the following assertion:
\begin{equation}\label{eqccs}
\vcenter{\it\hbox{Every complete algebraic system whose general member is
arithmetically}
 \hbox{effective and irreducible has a complete characteristic
system.}}
\end{equation}
Moreover, \eqref{eqccs} is equivalent to the existence of at least one
system with $R=r+q$, and in turn to the existence of suitably many such
systems.

On 16 January 1905 in the C. R. Paris, Enriques
\cite[\pp.134--135]{EnCR05} announced that Severi had just proved $q\ge
s$ and $q=b-s$, where $b$ is the number of independent Picard integrals
of the {\it second kind}, those with polar singularities.  It was known
before 1897 that $b$ is equal to the first Betti number; please see
\cite[\p.157]{Zas}.  In the same issue of the C. R., Picard \cite{Pcr05}
proved $q=b-s$ independently.%
 \footnote{Picard presented Enriques's note to the Academy, but
explained in Fn.\,$(^1)$ on \p.122 of \cite{Pcr05} that he had
completed his own note before  receiving Enriques's.}

 In the next issue on 23 January, Castelnuovo \cite{Ca1905} outlined the
last step in this direction.  He gave the details in three notes in the
Rend.\ Accad.\ Lincei of 21 May and 4 and 8 June 1905.  Specifically, he
took a complete algebraic system with arithmetically effective
 (in fact, regular)
 general member, fibered it into linear systems, and formed the quotient,
$P$ say.  Then $P$ is projective, and $P$ is of dimension $q$ as
$R=r+q$, an equation he considered proved.  Morover, $P$ is, up to
isomorphism, independent of the choice of algebraic system, and sum
(union) of curves induces an addition of points of $P$, turning $P$ into
a commutative group variety.

Hence, by a general 1895 theorem of Picard, completed in 1901 by
Painlev\'e, $P$ is an {\it Abelian variety}: $P$ is parameterized by $q$
{\it Abelian functions}, or $2q$-ply periodic functions of $q$
variables, with a common lattice of periods.  Castelnuovo proved that
these functions induce independent Picard integrals on $X$.  Therefore,
$q\le s$.  Thus Castelnuovo obtained the {\bf Fundamental Theorem of
Irregular Surfaces}:
\begin{equation*}\label{eqFTIS}
 	\dim P=q= s=b/2.
\end{equation*}
In 1905, the term ``Abelian variety'' was not yet in use.  So naturally
enough, Castelnuovo termed $P$ the {\it Picard variety\/} of $X$.%
 \footnote{Castelnuovo \cite[\p.221]{Ca1905} explained that ``out of
respect for Picard's profound research on surfaces [sic] admitting a
group of birational automorphisms, [he] proposes calling the variety $P$
 (and [a certain] group $G_d$) {\it the Picard variety (and Picard group)
associated to the surface $X$.}''

Andre Weil \cite[I, \p.572]{Wcp} discussed his own use of
the term ``Picard variety'' in his commentary on his 1950 paper on
Abelian varieties.  Weil said, ``Historically speaking, it would have
been justified to give it Castelnuovo's name, but it was a matter of
tampering as little as possible with common usage rather than rendering
due homage unto this master.''}

In 1903, Severi \cite[\S\,26]{Ssp} discovered a remarkable expression
for $p_a$ in terms of a different Hilbert polynomial.  Say the smooth
surface $\wt X$ is a subvariety of some higher dimensional projective
space $\mathbb P^N$.  Form the Hilbert polynomial $\wt F(\nu)$ of the
homogeneous ideal of $\wt X$.  Then $\wt F(0)-1=p_a$.

 Serre, in his 1954 ICM talk \cite[\pp.286--291]{ScpI}, announced a
theory of coherent algebraic sheaves, inspired by the analytic work of
Friedrich Hirzebruch, Kunihiko Kodaira, and Donald Spencer.  In
particular, Serre proved the Euler characteristic of the twisted
structure sheaf $\chi(\mc O_{\wt X}(\nu))$ is equal to $\smash{\wt
F(\nu)}$.  Thus $p_a=\chi(\mc O_{\wt X})-1$;%
 \footnote{Afterwards, it became common to  define $p_a$ by this
formula.}
  so $p_a$ is independent of
the embedding of ${\wt X}\subset\mathbb P^N$ and of the projection
${\wt X}\to\mathbb P^3$.

Further, Serre Duality yields this equality of dimensions of cohomology
groups: $\smash{h^i(\mc O_{\wt X})=h^{2-i}(\Omega^2_{\wt X})}$ for all
$i$ where $\smash{\Omega^2_{\wt X}}$ is the sheaf of algebraic
$2$-forms.  However, $p_g=\smash{h^0(\Omega^2_{\wt X})}$, essentially by
definition, and
 $\smash{h^0(\mc O_{\wt X})} =1$ as $\wt X$ is connected.  Thus
\begin{equation*}\label{eqSd}
p_a=\chi(\Omega^2_{\wt X})-1,\quad p_g=h^2(\mc O_{\wt X}),\quad
  q=h^1(\mc O_{\wt X}).
\end{equation*}

Above, the first equation is a form of Severi's discovery.  Here is a
proof \label{eqpa} of it using Grothen\-dieck's generalization of Serre
Duality.  Let $\omega_X$ be the dualizing sheaf.  Since $X\subset\mathbb
P^3$ and $\Omega^3_{\mathbb P^3}=\mc O_{\mathbb P^3}(-4)$, duality theory
and simple computation yield $$\omega_X= \textit{Ext}^1(\mc O_X\,,\
\Omega^3_{\mathbb P^3})=
        \mc O_X(n-4).$$
 Hence, duality for the finite map $\pi\:\wt X\to X$ and elementary
manipulation yield
$$\pi_*\Omega^2_{\wt X}=\sHom(\pi_*\mc O_{\wt X}\,,\,\,\omega_X)
        = \go C(n-4)\quad\text{where}\quad
        \go C := \sHom(\pi_*\mc O_{\wt X}\,,\,\mc O_X).$$
Thus $\chi(\Omega^2_{\wt X}) = \chi(\go C(n-4))$. 

Here, $\go C$ is the conductor; it is the ideal sheaf on $X$ of the
 curve $\Gamma$ of double points.  Let $\go C_0$ be the ideal sheaf on
$\mathbb P^3$ of $\Gamma$.  Then $p_a=\chi(\go C_0(n-4))$, essentially
by definition.  Form the standard exact sequence
$$0\to \mc O_{\mathbb P^3}(-4) \xrightarrow{\smash{\x f}}
         \go C_0(n-4) \to \go C(n-4) \to 0.$$
By Serre's Computation, $\chi(\mc O_{\mathbb P^3}(-4))=-1$.  Thus
$p_a =  \chi(\go C(n-4))-1$, as desired.  

Recall $q=h^1(\mc O_{\wt X})$.  Also, $s=h^0(\Omega^1_{\wt X})$
essentially by definition.  So Hodge Theory yields $q=s$ and $q=b/2$,
but Hodge Theory is not algebraic.  However, a $p$-adic algebraic proof
that $q=s$ was given by Kirti Joshi \cite{Jhs}.  Further, if by $b$ is
meant the dimension of the first \G\ \'etale cohomology group, then a
standard algebraic argument yields $\dim P =b/2$; see \cite[Lem.\,2A7,
\p.375]{KacWc} for example.  The latter argument also works in positive
characteristic, but the equations $\dim P=q$ and $q=s$ may fail.  In
1955, Jun-ichi Igusa gave an example with $\dim P=1$ but $q=s=2$; in
1958, Serre \cite[\p.529]{ScpI} gave one with $\dim P=s=0$ but
$q=1$.

Finally, $H^1(\mc O_{\wt X})$ is always the Zariski tangent space at $0$
to the Picard scheme by Grothendieck's theory, and over $\bb C$ the
Picard scheme is smooth by Cartier's theorem; so $\dim P=q$.  Thus the
Fundamental Theorem of Irregular Surfaces can be proved algebraically
over $\bb C$ , and the proof does not involve the Theorem of
Completeness of the Characteristic System.  Yet, the latter theorem has
taken on a life of its own, and \G's work is heavily involved in proving
both theorems algebraically.  All that work is discussed further in
Section 5.

\bigbreak
\section{The Picard Variety}
\begin{flushright}{\it
 Ever since 1949, I considered the construction\\
 of an algebraic theory of the Picard variety\\
 as the task of greatest urgency in abstract algebraic geometry.}\\
 Andr\'e Weil \cite[II, \p.537]{Wcp}
\end{flushright}
\nobreak
\medskip

Up to 1949, Weil worked primarily in Number Theory and Algebraic
Geometry.%
 \footnote{``Weil was far from confining himself to'' those subjects, as
Serre \cite[\pp.523--526]{Saw} noted, citing Weil's work in real
and complex analysis, representation theory, and differential geometry.}
   That work culminated in proofs of the Riemann hypothesis for curves
in 1948 and in the formulation of his celebrated conjectures for
arbitrary dimension in 1949.  Next, he led ``the construction of an
algebraic theory of the Picard variety.''  In turn, that theory led
Grothendieck to develop his theory of the Picard scheme.  However, the
Weil Conjectures themselves motivated much of \G's work.  In particular,
they led to the notion of a \G\ topology, which, as noted in the
introduction, is fundamental for the very definition of the Picard
functor; that functor is the subject of Section 4.

In the present section, so that we may better appreciate \G's advances,
let us consider in chronological order up to 1962, what was sought and
what was proved about the Weil conjectures and the Picard variety.  Good
secondary sources include Dieudonn\'e's history \cite{Dhdag} for all of
it, Mazur's 1974 expository article \cite{Mpspm29} for the Weil
conjectures, and the explanatory comments and historical notes in Serge
Lang's 1959 book \cite{Lab} for the Picard variety.  Again, as those
sources contain many primary references, those references are not always
repeated here.

In his 1921 thesis, which was published in 1924, Emil Artin developed an
analogue of the classical Riemann hypothesis, in effect, for a
hyperelliptic curve over a prime field of odd characteristic.  In 1929,
Friedrich (F-K) Schmidt generalized Artin's work to all curves over all
finite fields, recasting it in the algebro-geometric style of Richard
Dedekind and Heinrich Weber.  In 1882, they had viewed a curve as the
set of discrete valuation rings in a finitely generated field of
transcendence degree 1 over $\mathbb C$, but their approach works in any
characteristic.  In particular, Schmidt ported their proof of the
Riemann--Roch theorem, and used it to prove that Artin's Zeta Function
satisfies a natural functional equation.

In 1936, Helmut Hasse proved Artin's Riemann hypothesis in genus 1 via
an analogue over finite fields for the theory of elliptic functions.
Then he and Max Deuring noted that to extend the proof to higher genus
would require developing a similar analogue for the nineteenth century
theory of correspondences between complex curves.

  Their work inspired Weil.%
 \footnote{The interaction among the three and others has attracted a
lot of study.  One delightful and well-documented report was published
by Mich{\`e}le Audin in 2012 as \cite{Aulag}.  It describes the
political, social, and personal circumstances at the time, while
focusing on three reviews of Weil's first note.}
 In each of two notes, \cite[I, \pp.257--259]{Wcp} of 1940 and \cite[I,
\pp.277--279]{Wcp} of 1941, he sketched a different proof of the Riemann
hypothesis in any genus.  In both, the key is a certain positivity
theorem for correspondences.  It was found over $\mathbb{C}$ by
Castelnuovo%
 \footnote{In both notes, Weil cites only Severi.  In his commentary on
the second note \cite[I, \p.553]{Wcp}, Weil wrote, ``it's one of
Castelnuovo's most beautiful discoveries (see his {\it Memorie Scelete,}
no.\,XXVIII, \pp.509--517).  But I didn't read Castelnuovo until 1945 in
Brazil; then I realized that Severi in the {\it Trattato}
(\cite[\pp.286--287]{Str}) had not given his elder due credit.''}
  in 1906, and proved over a field of any characteristic by Weil in two
ways: in 1940 by algebraizing Adolf Hurwitz's transcendental theory of
1886, and in 1941 by porting to positive characteristic the
algebro-geometric theory in Severi's textbook \cite{Str} of 1926.

To provide the details, Weil had to redo the foundations of Algebraic
Geometry over a field of arbitrary characteristic.  The first instalment
\cite{Wf} appeared in 1946.  Building on ideas of Emmy Noether, Bartel
van der Waerden, and Schmidt from the 1920s, Weil fixed a {\it universal
domain} $\Omega$, a field of infinite transcendence degree over the
prime field.  Then a {\it projective variety\/} $X$ is the locus of
zeros with coordinates in $\Omega$ of homogeneous polynomials with
coefficients in a variable {\it coefficient field} or {\it field of
definition} $k$, a subfield of $\Omega$ over which $\Omega$ has infinite
transcendence degree.  Also $X$ is {\it absolutely irreducible}, not the
union of two smaller such loci.  Then Weil formed {\it abstract\/}
varieties by patching pieces of projective varieties.

Finally, Weil treated {\it cycles}.  They are the formal $\mathbb
Z$-linear combinations of subvarieties.  Those of codimension 1 are
called (\!{\it Weil\/}) {\it divisors}, and play a major role in the
theory of the Picard variety.  Weil developed a calculus of cycles,
including intersection products, inverse images, and direct images.

In 1948, Weil published two books.  In the first \cite{Wc}, he completed
his note of 1941.  He reproved the Riemann--Roch theorem, and developed
an elementary theory of correspondences for curves.  To prove
Castelnuovo's theorem, he used his full calculus of cycles on products
of numerous copies of the curve.  The proof is ``the most complicated
part of the'' book, as Otto Schilling observed in his Math Review
[MR0027151].  Then Weil proved the Riemann hypothesis.%
 \footnote{It is extraordinarily important.  Dieudonn\'e
\cite[\p.83]{Dhdag} gave one reason why: it ``allows proofs, in analytic
[sic] number theory, of `the best possible' upper bounds, 
inaccessible'' otherwise, such as this bound on a Kloosterman sum:
 $\bigl|\sum_{x=1}^{p-1}\exp\bigl(\frac{2\pi i}p(x+x^{-1})\bigr)\bigr|
  \le  p^{1/2}$ for any prime $p$.}

Weil's proof inspired three others.  First, in his 1953 thesis under
Hasse, Peter Roquette translated it into the more arithmetic language of
Schmidt, and simplified it to involve the product of just two
different curves.  Second, in 1958, Arthur Mattuck and John Tate applied
the Riemann--Roch theorem for surfaces, which had been proved in any
characteristic by Zariski in 1952 and by Serre in 1956.  Mattuck and
Tate proved the version of Castelnuovo's theorem for the product of two
curves that Severi \cite{Str} gave on \p.265.  They dubbed it the {\it
inequality of Castelnuovo--Severi}.  Then they rederived the Riemann
hypothesis, thus showing that it is a fairly simple consequence of the
general theory of algebraic surfaces.

Third, right as Mattuck and Tate finished their paper, \G\
 \cite[\p.208]{Gnmt}, ``attempting to understand the full import of
their method,'' found that it produces an index theorem on any surface,
which yields the Castelnuovo--Severi inequality.  According to \G\
however, Serre pointed out to him that he had proved an algebraic version
of William Hodge's 1937 analytic index theorem, and moreover that the
same version had already been proved the same way by Beniamino
Segre in 1937 and independently by Jacob Bronowski in 1938.

In Weil's second book \cite{Wva} of 1948, he completed his note of 1940.
He developed the abstract theory of {\it Abelian varieties}, which he
defined as the group varieties that are {\it complete}, the abstract
equivalent of ``compact.''  He proved that they are commutative, and
that a map between two is a homomorphism plus a translation.

Weil constructed the Jacobian $J$ of a smooth curve $C$ of genus $g$ by
patching together copies of an open subset of the symmetric product
$C^{(g)}$.  Given a prime $l$ different from the characteristic, he
constructed, out of the points on $J$ of order $l^n$ for all $n\ge1$, an
$l$-adic representation of the ring of correspondences, which, over
$\bb C$, is equivalent to the representation on the first
cohomology group.  He proved that the trace of this representation is
positive definite, and recovered Castelnuovo's theorem.  Finally, he
reproved the Riemann hypothesis for curves.

Weil left open, as Lang \cite[\p.17]{Lab} noted, two important
questions: (i) Is $J$ defined over the given coefficient field of {$C$}?
(ii) Is every Abelian variety projective?  Both questions were soon
answered affirmatively: (i) by Wei-Liang Chow, who announced his answer
in 1949 but published it in 1954, and (ii) by Teruhisa Matsusaka in
1953.  In 1954, Weil gave a much simpler and more direct answer to (ii).

In 1956, in order to handle (i), Weil addressed the general question of
finding a smaller coefficient field, but only in the case where the
resulting field extension is finitely generated and separable.  In turn,
Weil's work inspired \G\ to develop his general Descent Theory, which he
then sketched in his December 1959 Bourbaki talk \cite[190]{FGA}.
Grothendieck said on \p.190-1 that he was also inspired by Cartier's
subsequent treatment \cite[\S4]{Csb164} of purely inseparable
extensions, but that ``due to the lack of the language of schemes, and
especially the lack of nilpotents, Cartier could not express the basic
commonality of the two cases.''

In 1949, Weil published his celebrated conjectures about the zeta
function of a variety of arbitrary dimension.  He did not involve a
hypothetical cohomology theory outright, but one is implicit.  Moreover,
one was credited to him explicitly in Serre's 1956 ``Mexico paper''
\cite[\p.502]{ScpI} and in \G's 1958 ICM talk \cite[\p.103]{Gicm58},
where the term ``Weil cohomology'' appears, likely for the first time.

\G\ then announced ``an approach [to Weil cohomology]\dots sug\-ges\-ted
to [him] by the connections between sheaf-theoretic cohomology and
cohomology of Galois groups on the one hand, and the classification of
unramified coverings of a variety on the other\dots, and by Serre's idea
that a `reasonable' algebraic principal fiber space\,\dots should become
locally trivial on some covering unramified over a given point.''  Thus,
on \p.104, he announced the birth of {\it Grothendieck topology}.

In 1950, Weil published a remarkably prescient note \cite[I,
\pp.437--440]{Wcp} on Abelian varieties.  For each  normal%
 \footnote {{\it Normal\/} means the singular locus has codimension at
least 2, and (Zariski's Main Theorem) a rational map is defined
everywhere if its graph projects finite-to-one onto $X$ (so
isomorphically).}
 projective variety $X$ of any dimension in any characteristic, he said
that there ought to be two associated Abelian varieties, the {\it Picard
variety} $P$ and the {\it Albanese%
 \footnote {In 1913, Severi introduced and studied $A$ over the complex
numbers.  Nevertheless, much to Severi's dismay, Weil \cite[I,
\p.571]{Wcp} named $A$ after Severi's student, Giacomo Albanese,
ostensibly because, in 1934, Albanese viewed $A$ as a quotient of a
symmetric power of $X$.  However, Weil \cite[I, \p.562]{Wcp} left the
impression that rather it is because he owed a debt of gratitude to
Albanese for enriching the library of the University of Sa\~o Paulo,
Brazil, with works of Castelnuovo, Torelli and others, which were new to
Weil and from which he ``profited amply.''%
 }
  variety\/} $A$, with  the following  properties:
\begin{description}
 \item [Universality] The Picard variety $P$ parameterizes the linear
equivalence classes of all divisors on $X$ algebraically equivalent%
 \footnote{{\it Algebraic equivalence} and {\it  linear equivalence\/}
are just the equivalence relations generated by the algebraic  systems and
the linear systems.}
  to
0.  There is a  rational map%
 \footnote{A {\it rational map\/} is given by the ratio of two
polynomials, and is {\it defined\/} at a point if, for some choice of
the two,
the denominator does not vanishes there.}
 from $X$ into $A$, defined wherever $X$ is smooth, such that every
rational map from $X$ into an Abelian variety factors uniquely, up to
translation, through it.

\item [Duality] If $X$ is an Abelian variety, so that $X=A$, then $A$ is
the Picard variety of $P$; such a pair, $A$ and $P$, are called {\it
dual\/} Abelian varieties.
\end{description}
 Also, $A$ and $P$ are {\it isogenous}, or finite coverings of each
other, and of dimension equal to the irregularity [sic].  If $X$ is
arbitrary, then $A$ and $P$ are dual; in fact, the universal map $X\to
A$ induces the canonical isomorphism from the Picard variety of $A$ to
$P$.  If $X$ is a curve, then both $A$ and $P$ coincide with the
Jacobian.

In the note, Weil said that he had complete treatments of $P$ and $A$
for a smooth complex $X$, and sketches in general.  The sketches rest on
two criteria for linear equivalence of divisors in terms of linear-space
sections.  The criteria were found in 1906 by Severi, and reformulated
in the note by Weil, who referred to \pp.104--105,\, 164--165 in
Zariski's book \cite{Zas}; please see Mumford's comments
\cite[\p.120]{Zas} as well.  Weil announced proofs of the criteria in
1952, and gave the details in 1954.

In 1951, Matsusaka gave the first algebraic construction of $P$.  He
extend the coefficient field $k$ in order to apply two of Weil's
results: one of the equivalence criteria and the construction of the
Jacobian.  Both applications involve the {\it generic curve}, the
section of $X$ by a linear space of appropriate dimension defined over a
transcendental extension of $k$.  In 1952, Matsusaka gave a second
construction; it does not require extending $k$, but does require $X$ to
be smooth.

Both of Matsusaka's constructions are like Castelnuovo's: first Matsusaka
constructed a complete algebraic system of sufficiently
positive divisors, and then he formed the quotient modulo linear
equivalence.  To parameterize the divisors, he used the theory of
``Chow coordinates,'' which was developed in 1938 by Chow and van der
Waerden and was under refinement by Chow.  In 1952, Matsusaka also used
Chow coordinates to form the quotient.  Further, he made the first
construction of $A$, again using the Jacobian of the generic curve, but
he did not relate $A$ and $P$.

Also in 1952, in \S\,II of his paper {\it On Picard Varieties\/}
\cite[II, \pp.73--102]{Wcp}, Weil refined the sense in which $P$
parameterizes classes of divisors.  Working complex analytically, he
constructed ``an algebraic family of divisors on $X$, parameterized by
$P$, containing one and only one representative of each class.''

Weil did not name that family of divisors.
  However, the same year, Andr\'e N\'eron and
Pierre Samuel \cite{NSvpvn} constructed,%
 \footnote{Unfortunately, Lang \cite[\p.175]{Lab} felt that he had to write:
``It can not be said that the Picard variety is constructed in
\cite{NSvpvn} because this paper begins by a false statement concerning
the birational invariance.  This is a delicate point, when the varieties
involved have singularities.''}
 in any characteristic, a similar family, which they named a {\it
Poincar\'e family\/} citing \cite[II, \pp.73--102]{Wcp} in a way 
suggesting  the name%
 \footnote{Of course,  here and below, Poincar\'e's name  is used to honor
his work  mentioned above.}
  is due to Weil.  The family is defined by a divisor $D$ on $X\x P$,
which is called a {\it Poincar\'e divisor\/} by Lang \cite[\p.114]{Lab}.
Moreover, Lang showed that the pair $(P,D)$ is unique, $P$ up to
isomorphism and $D$ up to addition of a ``trivial'' divisor.

In 1955, Chow constructed $A$ and $P$ in a new way, as what he called
respectively the ``image'' and the ``trace'' of the Jacobian of a
generic curve on $X$.  Also, he proved that, indeed, the universal map
$X\to A$ induces an isomorphism of the Picard variety of $A$ onto $P$.

In a course at the University of Chicago, 1954--55, Weil gave a more
complete and elegant treatment, based on the ``see-saw principle,''
which he adapted from Severi, and on his own Theorem of the Square and
Theorem of the cube.  This treatment became the core of Lang's 1959 book
\cite{Lab}.  The idea is to construct $A$ first using the generic curve,
and then to construct $P$ as a quotient of $A$ modulo a finite subgroup.
Thus there is no need for Chow coordinates.

In 1958, Serre \cite[\p.555]{ScpI} worked over a fixed algebraically
closed coefficient field  $k$ of any characteristic.  He reproved
Igusa's 1955 bound $\dim A \le h^0(\Omega^1_X)$, and obtained a simple
direct construction of $A$ over $k$, not using the generic curve.

In 1958 Cartier \cite{Csb164} and in 1959 Nishi \cite{Nd} independently
proved Weil's conjectured duality of $A$ and $P$: in any
characteristic, each is the others Picard  variety.

Between 1952 and 1957, Maxwell Rosenlicht published a remarkable series
of papers, inspired by Severi's 1947 monograph \cite{Sfqa}, which
treated curves with double points.  Treating a curve with arbitrary
singularities, Rosenlicht generalized the notions of linear equivalence
and differentials of the first kind.  Then he constructed a {\it
generalized Jacobian} $J$ over $\mathbb{C}$ by integrating and in
arbitrary characteristic by patching.  It is not an Abelian variety, but
an extension of the Jacobian $J_0$ of the desingularized curve by an
affine algebraic group.  In 1962, Frans Oort \cite{Oge} gave another
construction, which gives $J$ as a successive extension of $J_0$ by
additive and multiplicative groups.

For arithmetic applications, Tate suggested, according to Lang
\cite[\p.176]{Lab}, doing this.  Given finitely many simple points on
$X$, consider the divisors avoiding them.  Form linear equivalence
classes via functions congruent to 1 to given multiplicities at the
points.  Finally, seek a generalized Picard variety to parameterize
these classes.

In 1959, Serre published a textbook \cite{Sgacc} on the case $\dim X =1$
and its application to Lang's class field theory over function fields.
In particular, Serre recovered Rosenlicht's generalized Jacobian of an
$X$ with one singular point%
 \footnote{The singularity cannot be arbitrary; for example, it cannot
be a planar triple point.}%
 , constructed by identifying given points
with given multiplicities on a given smooth curve.

In 1962, Murre \cite{MgPv} constructed Tate's generalized Picard variety
$P$ by adapting Matsusaka's two constructions.  Thus Murre obtained $P$
for any (normal) $X$ via patching and for any smooth $X$ directly over
the same ground field.

In 1956, Igusa established the compatibility of specializing a curve
with specializing its generalized Jacobian, possibly under reduction mod
$p$, provided the general curve is smooth and the special curve has at
most one node.  Igusa explained that, in 1952, N\'eron had studied the
total space of such a family of Jacobians, but had not explicitly
analyzed the special fiber.%
 \footnote{For a comprehensive discussion of the ``N\'eron model'' and its
connection to the Picard scheme (and Picard algebraic space) along with
historical notes and references to the original sources, please see the
textbook \cite{BLRnm} of Siegfried Bosch, Werner L{\"u}tkebohmert, and
Raynaud.}
   Igusa's approach is, in spirit, like Castelnuovo's, Chow's, and
Matsusaka's before him.

In 1960, Claude Chevalley \cite{Ctp} constructed a Picard variety using
Weil divisors locally defined by one equation; they are called {\it
Cartier divisors} in honor of Cartier's 1958 Paris thesis \cite{Cbsm}.
First, Chevalley constructed a {\it strict\/} Albanese variety; it is
universal for {\it regular\/} maps, ones defined everywhere, into
Abelian varieties.  Then he took its Picard variety to be that of $X$.
He noted his Picard and Albanese varieties need not be equal to those of
a desingularization of $X$.  By contrast, Weil's $P$ and $A$ are
birational invariants, and his universal map $X\to A$ is a rational map.
In 1962, Conjeerveram Seshadri \cite{SvP} generalized Chevalley's
construction to an $X$ with arbitrary singularities, 
recovering Rosenlicht's generalized Jacobian.

In 1961, Mattuck \cite{Mpb} \label{Mat1} took, on a smooth $X$ over an
algebraically closed field, a complete algebraic system $\Sigma$ of
suitably positive divisors $C$.  He parameterized $\Sigma$ by the {\it
Chow variety\/} $H$, the locus of points given by the Chow coordinates
of the $C\in \Sigma$.  He fixed a $D\in\Sigma$, and took the class of
the difference $C-D$ for $C\in\Sigma$, to get a rational map
$\alpha\:H\to P$.  In order that $\alpha$ be defined everywhere, he
reembedded $X$ in another projective space, because Murre \cite{MCv} had
just proved that then $H$ is smooth, so normal.

  Mattuck proved that $\alpha$ is a projective bundle%
 \footnote{Earlier, in 1956, Kodaira \cite{Kcls} obtained a similar
result analytically.}
 and has a section.  The section corresponds to a refined Poincar\'e
divisor: not only does it define a Poincar\'e family, but it contains no
fiber of $\alpha$, so cuts each fiber in a divisor.  Finally, Mattuck
studied the case that $X$ is a curve of genus $g>1$, so that $H$ is the
$n$-th symmetric power of $X$ where $n$ is the degree of the divisors.
He proved that $\alpha$ is a bundle%
 \footnote{Atiyah \cite[\p.451]{Avbec} wrote in 1957 that this case is
``well known.''}
  if $n>2g-2$, and that $\alpha$ has a section%
 \footnote{In fact, $n>2g-1$ suffices; please see \p.\pageref{Mat}
below.}
 if $n>4g$, but no section if $n=2g-1$ and the divisor classes modulo
algebraic equivalence on $P$ form a group of $\mathbb{Q}$-rank 1.  So it
seems unlikely that $\alpha$ has a section if $n=2g-1$ and $X$ has
general moduli.

Thus, in 1962, the algebraic theory of the Picard variety was indirect,
involved, and incomplete.  There were competing definitions and
constructions, each with advantages and disadvantages.  There was a lot
of fussing with fields of definition.  There were loose ends.  Notably,
there was no fully satisfactory way to parameterize divisors or to construct
quotients.  So there was not enough machinery to prove the Completeness
of the Characteristic System or to treat, in general, the behavior
of the (generalized) Picard variety in a family.
\G\ brilliantly handled those issues in the way explained in the next
two sections.

\bigbreak
\section{The Picard Functor}
\begin{flushright}{\it
 Grothendieck certainly did not feel\\
 that he was attempting to use powerful techniques\\
 in order to obtain stronger results by generalizing.\\
 What he perceived himself as doing was simplifying situations and objects,\\
 by extracting the fundamental essence of their structure.}\\
 Leila Schneps \cite[\p.3]{SchnepsCh1}
\end{flushright}
\nobreak\medskip

Many times, \G\ proceeded ``by extracting the fundamental essence'' of
existing theories, and then developing his own versions, for example his
theories of schemes, representable functors, the Hilbert scheme, and the
Picard scheme.  Parts of those theories belong to the theory of
the Picard functor.  Those parts are treated in depth in \cite{EGAIG}
and \cite{FAG}, and they are introduced in this section.

To begin, here is a bit more informal history.  Starting in 1937,
Zariski made deep use of the local ring of all rational functions that
are finite at a given point of a variety with a fixed field of
definition.  Inspired by Zariski's work, Chevalley developed an
intrinsic theory of abstract varieties $V$ in his paper \cite{Cva},
submitted on 2 July 1954.  On \p.2, he called the set of all the local
rings of the points of $V$ its {\it model,} and then developed a theory
of models.  Earlier, in January 1954, he lectured on his theory at Kyoto
University, according to Masayoshi Nagata \cite[\p.78]{Najm56}, who then
generalized it by replacing the field of definition by a Dedekind
domain.

In 1944, Zariski topologized the set of all valuation rings of the field
of rational functions of a variety in order to use the finite-covering
property to pass from local uniformization to global desingularization.
In 1949, Weil \cite[I, \pp.411--413]{Wcp} observed that his abstract
varieties support what he called the {\it Zariski Topology}, whose
closed sets are the subvarieties and their finite unions.

Weil used the Zariski topology to define locally trivial fiber spaces.
He discussed the natural bijective correspondence between the line
bundles on a smooth variety and its linear equivalence classes of
divisors.  Then, in his 1950 paper on Abelian varieties \cite[I,
\pp.438--439]{Wcp}, he suggested that those line bundles might be used to
develop, for any abstract variety, a version of Severi's generalized
Jacobian.

As already noted in Section 2, Serre, in his 1954 ICM talk, announced a
theory of coherent algebraic sheaves.  In fact, he worked over an
arbitrary algebraically closed field $k$, of any transcendence degree
over the prime field, and used $k$ both as the field of definition and
as the field of coordinates.%
 \footnote{Cartier  \cite{Cbsm} generalized Serre's theory to an
arbitrary  field of definition, and studied the effect of extending it.
He took the  field of coordinates to be a universal domain.}
   Moreover, he worked only with projective space and its subvarieties,
which he allowed to be reducible, and he viewed as the closed sets of a
topology, which he too called the ``Zariski topology.''

In 1955, Serre presented the details in his celebrated paper {\it
Faisceaux alg\'ebriques coh\'erents\/} \cite[\pp.310--391]{ScpI}.  He
also generalized the notion of abstract variety via Henri Cartan's
notion%
 \footnote{Cartan had used the notion to define  $C^\infty$-manifolds.
He and Serre had used it to define complex analytic manifolds.}
  of {\it ringed space}.  It is a topological space $X$ endowed with a
{\it sheaf of rings} $\mc O_X$, called the {\it structure sheaf\/}: over
each open set, its sections form a ring; for each smaller open set, the
restriction map is a ring homomorphism.  To be a variety, $X$ must be
covered by finitely many open subsets, each of which, when endowed with
the restriction of $\mc O_X$, is isomorphic to an {\it affine
variety\/}; the latter's space is a closed subspace of an affine space,
and its structure sheaf has, as its sections over an open set $U$, the
rational functions defined everywhere on $U$.

Both Serre and Chevalley downgraded rational maps, preferring maps%
 \footnote{Serre called them {\it regular maps,} a traditional term;
Chevalley \cite[\p.219]{Csb58} used {\it morphisms}.}
 that are defined everywhere.  For Serre, a map of varieties $\vf\:X\to
Y$ is a map of ringed spaces: $\vf$ is a continuous map equipped with a
map $\vf^*$ relating the two structure sheaves, so that a section $f$ of
$\mc O_Y$ over an open set $V$ yields a section $\vf^*f$ of $\mc O_X$
over $f^{-1}V$ in a way respecting addition, multiplication, and
restriction.

\G\ ``extracted the fundamental essence of''\, those ideas, and
developed a theory of {\it  schemes.}%
 \footnote{Cartier \cite[Fn.\,8]{Ch} noted, ``This word results from a
typical epistemological shift from one thing to another: for Chevalley,
who invented the name in 1955, it indicated the `scheme' or `skeleton'
of an algebraic variety, which itself remained the central object.  For
\G, the `scheme' is the focal point, the source of all the projections
and all the incarnations.''}
 By February 1956 (see \cite[\p.32]{CGS}), he was working with ringed
spaces that have an open covering by {\it affine schemes}, or {\it ring
spectra}.%
 \footnote{Cartier \cite[Fn.\,29]{Ch} explained, ``It was [Israel]
Gelfand's fundamental idea [of 1938] to associate a normed
commutative algebra to a space.\,\ldots The term `spectrum' comes
directly from Gelfand.''}
   The spectrum of a ring $R$ is the set of all its prime ideals $\go
p$.  Its topology is generated by its {\it principal open subsets}
$D(f)$ for all $f\in R$, where $D(f):=\{\go p\ni f\}$.  Over $D(f)$, the
sections of the structure sheaf are the fractions $a/f^n$ for all $a\in
R$ and $n\ge0$.  In 1956, \G\ took $R$ to be Noetherian, but in EGA\,I
\cite{EGAI}, which appeared in 1961, $R$ is an arbitrary commutative
ring.  In any case, $R$ may have nilpotents.

Moreover, \G\ worked with Cartier's generalization in
\cite[\p.206]{Cbsm} of coherent sheaves, the {\it quasi-coherent}
sheaves.  He \cite[\p.32]{CGS} told Serre why: they are ``technically
very convenient because they have the relevant properties of coherent
sheaves without requiring the finiteness (on an affine, they correspond
to all the modules over the coordinate ring, and not just the finitely
generated modules).''

The generality is vast, but not idle.  Murre put it as follows,
according to Schneps \cite[\p.2]{SchnepsCh1}: ``Un\-doubt\-edly, people
did see in the mid 50's that one could generalize a lot of things to
schemes, but \G\ saw that such a generalization was not only possible
and natural, but necessary to explain what was going on, even if one
started with varieties.''

The spectrum $S$ of a field $k$ is a single point, but $S$ has $k$ as
structure sheaf.  Call $k$ the {\it field of definition\/} of a scheme
$X$ if there is a distinguished map $X\to S$.  On the other hand, given
a universal domain $\Omega$ extending $k$, let $T$ be its spectrum.
Then a {\it point\/} of $X$ with ``coordinates'' in $\Omega$, or an
$\Omega$-{\it point\/} of $X$, is just a map $T\to X$ that respects the
distinguished maps $X\to S$ and $T\to S$.  Cartier \cite[\p.20]{Ch}
described this new situation as being ``admirable simplicity\emdash and
a very fruitful point of view\emdash but a complete change of
paradigm!''

For example, say $X$ is the spectrum of the residue ring $R$ of the
polynomial ring $k[u_1,\dotsc,u_n]$ modulo the ideal generated by
polynomials $f_\lambda$.  Then the inclusion $k\to R$ defines the map
$X\to S$.  Further, since $\Omega$ is a universal domain, a zero
$(c_1,\dotsc,c_n)$ of the $f_\lambda$ in $\Omega^n$ amounts exactly to a
$k$-algebra map $\gamma\:R\to\Omega$ that carries the residue of $u_i$ to
$c_i$. In turn, $\gamma$ amounts exactly to an  $S$-map $T\to X$.

However, there's no need to restrict $S$ and $T$ to be the spectra of
fields, and good reason not to.  Moreover, the right setting for this
theory, as for any theory whose principals are objects and maps, is
Category Theory.  It is not simply a convenient language for expressing
abstract ideas, but more importantly, an effective tool, which eases the
work at hand and affords new possibilities.  \G\ recognized as much, and
promoted Category Theory.

Thus, given a {\it base\/} scheme $S$, an $S$-{\it scheme\/} is a scheme
equipped with a  map to $S$, its {\it structure map}.  An
$S$-{\it map} is a map between $S$-schemes that commutes with the two
structure maps.  The category of $S$-schemes has products:%
 \footnote{Mac Lane \cite[\p.76]{Mcwm} said that he himself, in 1948 and
1950, formulated ``the idea that\,\ldots\break products could be described by
universal properties of their projections.''}
  the product of $X$ and $Y$ is an $S$-scheme $X\x_SY$ with a
distinguished pair of  $S$-maps to $X$ and $Y$, the {\it projections},
 such that composition with them sets up a bijection from the $S$-maps
$T\to X\x_SY$ to the pairs of $S$-maps $T\to X$ and $T\to Y$.

The product $X\x_SY$ is determined, formally, up to unique
isomorphism.  It is constructed by patching together the spectra of the
tensor products of the rings of affines that cover $S$, $X$, and $Y$.
It can also be viewed as the result of {\it base change} of $X$ when $Y$
is viewed as a new base.  By contrast, Cartier \cite[\p.20]{Ch} noted
that, ``in both [Serre's and Chevalley's] cases, the two fundamental
problems of products and base change could only be approached
indirectly.''

Another way to view a map $X\to S$ is as a family $X/S$ with $S$ as
parameter space and $X$ as total space.  Its members are the {\it
fibers}, the preimages $X_s$ of the points $s\in S$.  More precisely,
$X_s:=X\x_SY$ where $Y$ is the spectrum of the residue field $k_s$ of
the stalk $\mc O_{S,s}$, which is  the local ring of all functions that are
each defined on some neighborhood of $s$.
  \G\ discovered that, many properties of the $X_s$ vary continuously
when $X$ is $S$-{\it flat}; that is, $\mc O_{X,x}$ is $\mc O_{S,s}$-flat
for all $s\in S$ and $x\in X_s$.  Often, he considered what he called a
{\it geometric fiber}, which is the product $X\x_S Y'$ where $Y'$ is the
spectrum of an algebraically closed field containing $k_s$.

An $S$-map $T\to X$ is called a $T$-{\it point\/} of $X$, and the set of
all of them is denoted by $X(T)$ or $h_X(T)$.  An $S$-map $T'\to T$
induces a set map $h_X(T)\to h_X(T')$.  Thus $h_X$ is a contravariant
functor from the category of $S$-schemes to the category of sets; it is
called the {\it functor of points\/} of $X$.

The contravariant functors $H$ from $S$-schemes to sets themselves form
a category.  The assignment $X\mapsto h_X$ is a functor from $S$-schemes
into the latter category.  This functor is an embedding by Yoneda's
Lemma.  Given an $H$, if an $X$ is found with $h_X=H$, then $H$ is said
to be {\it representable\/} by $X$.  If so, then $H(X)$ contains a {\it
universal\/} element $W$, which corresponds to the identity map of $X$;
namely, each element $Y\in H(T)$ defines a unique $S$-map $\vf\:T\to X$
with $H(\vf)W=Y$.
  In other words, the $T$-points of $X$ {\it represent\/} the elements
of $H(T)$.  Conversely, if an $X$ is found with a universal $W\in H(X)$,
then there's a canonical isomorphism $h_X=H$.

The first important example is the functor $P(\mc E)$, where $\mc E$ is
an arbitrary quasi-coherent sheaf on $S$.  For each $S$-scheme $T$, the
set $P(\mc E)(T)$ is the set of invertible quotients $\mc L$ of the
pullback $\mc E_T$; {\it invertible} means that $\mc L$ is the sheaf of
sections of a line bundle.  The functor $P(\mc E)$ is representable by
an $S$-scheme $\I P(\mc E)$.  Automatically, $\I P(\mc E)$ carries a
{\it universal\/} invertible quotient of $\mc E_{\I P(\mc E)}$, denoted
$\mc O_{\I P(\mc E)}(1)$; namely, each invertible quotient $\mc L$ of
$\mc E_T$ defines a unique $S$-map $\vf\:T\to \I P(\mc E)$ with
$\vf^*\mc O_{\I P(\mc E)}(1)=\mc L$.  Moreover, $\I P(\mc E)\x_SY=\I
P(\mc E_{Y})$ for any $S$-scheme $Y$.  In particular, the fiber $\I
P(\mc E)_s$ over $s\in S$ is the projective space of $1$-dimensional
quotients of the vector space $\mc E_s\ox k_s$ where $\mc E_s$ is the
stalk of $\mc E$  at $s$.

For convenience, {\bf assume from now on} that all schemes are {\it
locally Noetherian}, or covered by affine schemes on Noetherian rings,
and that $S$ is {\it Noetherian}, or covered by finitely many such.
Assume {\bf also} that $X$ is projective over $S$; that is, $X$ can be
embedded as a closed subscheme of $\I P(\mc E)$ for some coherent sheaf
$\mc E$ on $S$.

Consider the {\it Hilbert functor\/} $\Hilb_{X/S}$, treated by \G\ in
his May 1961 Bourbaki talk \cite[Exp.\,221]{FGA}.  For each $S$-scheme
$T$, the set $\Hilb_{X/S}(T)$ is the set of $T$-flat closed subschemes
$Y$ of $X\x_ST$.  Moreover, for each polynomial $F(\nu)$ with rational
coefficients, $\Hilb_{X/S}$ has a subfunctor $\Hilb_{X/S}^F$; namely,
for all $T$, the set $\smash{\Hilb_{X/S}^F(T)}$ is the set of $Y$ such
that $Y_t$ has Hilbert polynomial $F$ for all $t\in T$, that is,
$F(\nu)=\chi(\mc O_{Y_t}(\nu))$ where $\mc O_{Y_t}(\nu)$ is
the pullback to $Y_t$ of the $\nu$th tensor power of $\mc O_{\I P(\mc
E)}(1)$.  \G\ proved that $\Hilb_{X/S}$ is representable by a locally
Noetherian $S$-scheme $\hilb_{X/S}$, the {\it Hilbert scheme}.  In fact,
$\hilb_{X/S}$ is the disjoint union of projective $S$-schemes
$\smash{\hilb_{X/S}^F}$, which represent the functors $\smash{\Hilb_{X/S}^F}$.

Automatically, $X\x_S\hilb_{X/S}$ has a {\it universal\/} subscheme $W$;
namely, each $T$-flat closed subscheme $Y$ of $X\x_ST$ defines a unique
$S$-map $T\to \hilb_{X/S}$ with $W\x_{\hilb_{X/S}}T=Y$.  Note that the
$\smash{\hilb_{X/S}^F}$ depend on the choice of embedding of $X$ in some
$P(\mc E)$, but $\hilb_{X/S}$ does not.  Thus the Hilbert scheme is a
noble replacement for Chow coordinates; the latter only parameterize the
cycles on a variety $V$, and depend on the choice of embedding of
$V$ in projective space.

A subscheme $R$ of $X\x_SX$ defines a {\it flat and projective
equivalence relation\/} if each projection $R\to X$ is flat and
projective and if, for each $S$-scheme $T$, the subset $h_R(T)$ of
$h_X(T)\x h_X(T)$ defines a set-theoretic equivalence relation.
\G\ found two constructions of a {\it quotient\/} $X/R$ in the strongest
sense of the term.  Namely, first, an $S$-map $X\to Z$ factors through
$X/R$ if and only if the two compositions $R\rra X\to Z$ are equal; if
so, then $X/R\to Z$ is unique.  So by ``abstract nonsense,'' $X/R$ is
determined up to unique isomorphism.  Second, the quotient map $X\to
X/R$ is flat and projective, and the canonical map $R\to X\x_{X/R}X$ is
an isomorphism.%
 \footnote{It follows that $h_{X/R}$ is the fpqc sheaf associated to
$h_X/h_R$ in
the sense discussed below.}

\G's first construction \cite[\p.212-15]{FGA} uses quasi-sections to
reduce to the case where $X$ is affine and each $R\to X$ has finite
fibers.  However, his second construction \cite[\p.232-13]{FGA} is
easier and more elegant.  It proceeds as follows: $R$ lies in
$\Hilb_{X/S}(X)$, so defines a map $\vf\:X\to \hilb_{X/S}$; the graph
$\Gamma_\vf$ is a closed subscheme of the universal subscheme $W$;
finally, by \G's Descent Theory, $\Gamma_\vf$ descends to a closed
subscheme of $\hilb_{X/S}$, which is the desired $X/R$.

 {\bf Assume from now on} that $X$ is {\bf also} $S$-flat.  Then
$\Hilb_{X/S}$ has an important subfunctor $\Div_{X/S}$; namely, for each
$S$-scheme $T$, let $\Div_{X/S}(T)$ consist of the {\it effective
Cartier divisors} $D$ in $\Hilb_{X/S}(T)$, that is, the flat subschemes
$D$ whose ideal $\mc I_D$ is locally generated by one nonzerodivisor;
equivalently, $\mc I_D$ is invertible as an abstract sheaf. \G\
\cite[\p.232-10]{FGA} proved that $\Div_{X/S}$ is representable by an
open subscheme $\IDiv_{X/S}$ of $\hilb_{X/S}$.

Given an invertible sheaf $\mc L$ on $X$, define%
 \footnote{This subfunctor appears in \cite[\p.232-10]{FGA}, but the
notation for it comes from  \cite[\p.93]{Mcs}.} 
 a subfunctor $\LS_{\mc L/X/S}$ of $\Div_{X/S}$: for each $S$-scheme
$T$, let $\LS_{\mc L/X/S}(T)$ consist of the $D$ in $\Div_{X/S}(T)$ for
which there is an invertible sheaf $\mc M$ on $T$ such that the inverse
$\mc I_D^{-1}$  is isomorphic to the tensor product
on $X\x_ST$ of the pullbacks of $\mc L$ and $\mc M$.

 {\bf Assume in addition from now on} the geometric fibers of $X/S$ are
{\it integral}; that is, each affine ring of each geometric fiber is an
integral domain.  \G\ \cite[\p.232-11]{FGA} proved that then there is a
coherent sheaf $\mc Q$ on $S$, determined up to unique isomorphism, such
that $\LS_{\mc L/X/S}$ is representable by $\I P(\mc Q)$.  Hence, $\I
P(\mc Q)$ is equal to a closed subscheme of $\IDiv_{X/S}$.  Also, if
$H^1(\mc L|X_s)=0$ at $s\in S$, then $s$ has a neighborhood on
which $\mc Q$ is {\it free}, or isomorphic to $\mc O_S^r$ for some $r$.

In general, what makes a functor $H$ representable?  Say $H=h_X$.  Then
given any $S$-scheme $T$ and any open covering $\{T_\lambda\}$ of $T$,
two maps $T\to X$ are equal if their restrictions to each $T_\lambda$
are equal.  Furthermore, maps $\vf_\lambda\:T_\lambda\to X$ are the
restrictions of a single map $T\to X$ if, for all $\lambda$ and $\mu$,
the restrictions of $\vf_\lambda$ and $\vf_\mu$ to $T_\lambda\cap T_\mu$
are equal.  In other words, as $U$ ranges over the open sets of $T$, the
$H(U)$ form a sheaf.  The latter condition does not explicitly
involve $X$.  So it makes sense for any $H$, representable or not.  If it
is satisfied, $H$ is called a {\it Zariski sheaf}.

Here's another formulation.  Let $T'$ be the disjoint union of the
$T_\lambda$, and consider the induced map $T'\to T$.  Then $T'\x_TT'$ is
the disjoint union of the $T_\lambda\cap T_\mu$.  So the condition to be
a Zariski sheaf just means that the induced sequence of sets
\begin{equation}\label{eqess}
H(T) \to H(T') \rra H(T'\x_TT')
\end{equation}
is {\it exact\/}; that is, the first map is injective, and its image
consists of the elements of $H(T')$ whose two images are equal in
$H(T'\x_TT')$.

\G's Descent Theory yields more.  Let $T'\to T$ be an {\it fpqc\/} map;
namely, it is flat and surjective, and the preimage of any affine open
subscheme is a finite union of affine open subschemes.  If $H$ is
representable, remarkably \eqref{eqess} is still exact; in other words,
$H$ is an {\it fpqc sheaf}.  Indeed, the {\it fpqc \G\ topology\/} may
be defined as the refinement of the Zariski topology with the fpqc maps
as additional generalized open coverings.  In particular, $H$ is an {\it
\'etale sheaf}, the notion obtained by requiring the maps to be
{\it \'etale}, that is, flat, unramified and locally of finite type.

The {\it Picard group} $\Pic(X)$ is the group, under tensor product, of
isomorphism classes of invertible sheaves on $X$.  The {\it absolute
Picard functor} $\Pic_X$ is defined by $\Pic_X(T):= \Pic(X\x_ST)$.  It
is never a Zariski sheaf, so never representable.

There is a sequence of ever more promising ``Picard functors.''  First
comes the {\it relative Picard functor} $\Pic_{X/S}$ defined by
$$\Pic_{X/S}(T):=\Pic(X\x_ST)/\Pic(T)$$
 where $\Pic(T)$ acts via pullback.  Following it are its {\it
associated sheaves\/} in the Zariski, \'etale and fpqc topologies:
$\Pic_{(X/S)\zar}$, $\Pic_{(X/S)\et}$, and $\Pic_{(X/S)\fpqc}$.  \G\
\cite[\p.232-03, (1.6)]{FGA} formed them as direct limits.  For example,
 $$\textstyle \Pic_{(X/S)\zar}(T) := \varinjlim_{T'}\Pic_{X/S}(T')$$
 where $T'$ ranges over the small category of all open coverings of $T$.

 {\bf Recall} that $S$ is Noetherian and that $X$ is a flat and projective
$S$-scheme with integral geometric fibers.  \G\ \cite[\pp.232-4--6]{FGA}
proved%
 \footnote{In \cite{FGA}, \G\ did not consider $\Pic_{(X/S)\et}$.
However, his methods apply to it, and show that it is equal to
$\smash{\Pic_{(X/S)\fpqc}}$, because in the case at hand, there exists
an {\it \'etale quasi-section}, an $S$-map $S'\to X$ for which the
structure map $S'\to S$ is
\'etale.}
 that {\it then the three canonical comparison maps are, respectively,
injective, injective, and bijective:
$$\Pic_{X/S}\into\Pic_{(X/S)\zar}\into\Pic_{(X/S)\et}\risom\Pic_{(X/S)\fpqc}.$$
Moreover, the first two maps are bijective if $X\to S$ has a section;
the middle map is bijective if it just has a section on a Zariski
neighborhood of each point of $S$.}

A simple example shows that, in general, we must pass to the \'etale
sheaf.  Namely, in the real projective plane, consider the conic $X :
u^2+v^2+w^2=0$.  Let $S$ and $T$ be the spectra of $\bb R$ and $\bb C$.
Then $T\to S$ is an \'etale covering.  Moreover, $X\x_ST$ is the complex
conic with the same equation; so $X$ is isomorphic to the complex
projective line.  The latter's universal sheaf $\mc O(1)$ defines an
element $\tau\in\Pic_{(X/S)\et}(S)$, as the two pullbacks of $\mc O(1)$
to $X\x_ST\x_ST$ are isomorphic.  And $\tau$ is not in the image of
$\Pic_{(X/S)\zar}(S)$, as $X$ has no $S$-points. 

\G's main existence theorem \cite[\p.232-06]{FGA} says that
{\it$\Pic_{(X/S)\et}$ is representable by a scheme $\pic_{X/S}$.}  It is
called the {\it Picard scheme}.  Of course, if $\Pic_{(X/S)\zar}$ is
representable, then it is an \'etale sheaf, so equal to
$\Pic_{(X/S)\et}$, and representable by $\pic_{X/S}$.  Similarly, if
$\Pic_{X/S}$ is representable, then all four functors are equal, and
representable by $\pic_{X/S}$.

\G's proof is fairly simple at this point.  Here is the idea.  Fix an
embedding of $X$ in a $\I P(\mc E)$.  Given any $S$-scheme $T$ and any
quasi-coherent sheaf $\mc F$ on $\I P(\mc E)\x_ST$, let $\mc F(n)$
denote the tensor product of $\mc F$ and the pullback of the $n$th
tensor power of the universal sheaf $\mc O_{\I P(\mc E)}(1)$.  Let $\mc
I$ be the ideal of the universal divisor on $X\x_S\IDiv_{X/S}$, and $\mc
I^{-1}$ its inverse.  Form the open subscheme $\I D^+$ of
$\IDiv_{X/S}$ on which all the higher direct images of $\mc I^{-1}(n)$
vanish for all $n\ge0$.

Set $\mc L:=\mc I^{-1}|(X\x_S\I D^+)$.  Then $\LS_{\mc L/X\x_S\I D^+/\I
D^+}$ is representable by $\I P(\mc Q)$ where $\mc Q$ is a coherent
sheaf on $\I D^+$.  Moreover, $\mc Q$ is {\it locally free\/}; that is, each
point of $\I D^+$ has neighborhood on which $\mc Q$ is free.  Thus $\I
P(\mc Q)$ is flat over $\I D^+$.

Set $R:=\I P(\mc Q)$.  Then $R$ is a closed subscheme of $\I D^+\x_S\I
D^+$.  Moreover, for each $S$-scheme $T$, the subset $h_R(T)$ of $h_{\I
D^+}(T)\x h_{\I D^+}(T)$ consists of the pairs of $D,\,D'\in \I D^+(T)$
for which there is an invertible sheaf $\mc M$ on $T$ such that the
ideal $\mc I_D$ is isomorphic to the tensor product of $\mc I_{D'}$ and
the pullback of $\mc M$.  Thus $h_R(T)$ is a set-theoretic equivalence
relation.

Although $\I D^+$ isn't Noetherian, nevertheless it is the disjoint
union of Noetherian subschemes, as $\hilb_{X/S}$ is the disjoint union
of the projective $S$-schemes $\hilb_{X/S}^F$, and $R$ decomposes
compatibly.  Consequently, the quotient $\I D^+/R$ exists.

For each $m\ge0$, let $P_m$ be the fpqc subsheaf of $\Pic_{(X/S)\fpqc}$
associated to the subfunctor of $\Pic_{X/S}$ whose value at $T$ consists
of the classes of invertible sheaves $\mc L$ on $X\x_ST$ for which all
the higher direct images of $\mc L(n)$ vanish on $T$ for all $n\ge m$,
but the direct image doesn't vanish.  Then $P_0$ is representable by $\I
D^+/R$.

Tensor product with the pullback of $\mc O_X(1)$ defines an isomorphism
$P_{m+1}\risom P_m$ for all $m\ge0$.  So the $P_m$ are representable by
(isomorphic) schemes $U_m$.  Each inclusion $P_m\into P_{m+1}$ is
representable by an open embedding $U_m\into U_{m+1}$.  Finally,
$\Pic_{(X/S)\fpqc}$ is the ``union'' of the $P_m$; so is representable
by the union, or rather direct limit, of the $U_m$.  Thus \G\ proved his
main existence theorem.

Commenting on his proof, \G\ \cite[\p.232-13]{FGA} noted that ``the
approach is basically the one followed by Matsusaka'' (so by Igusa, Chow
and Castelnuovo).  Further, he \cite[\p.232-14]{FGA} noted that ``the
proof appeals neither to the preliminary construction of the Jacobians
of curves\,\ldots\,nor to the theory of Abelian varieties, and thus
differs in an essential way from the `traditional' treatments in Lang's
book \cite{Lab} and Chevalley's paper \cite{Ctp}.\,\ldots That the
construction of the Picard scheme ought to precede and not follow the
theory of Abelian varieties is clear a priori from the fact that\,\dots
Rosenlicht's `generalized Jacobians' are not Abelian varieties.''  More
of \G's advances are highlighted in the next section.

\bigbreak
\section{The Picard Scheme}
\begin{flushright}{\it
 His [\G's] feeling was that ``those people''\\
 made too strict assumptions and tried to prove too little.}\\
 Jacob Murre, quoted in \cite[\p.2]{SchnepsCh1}
\end{flushright}
\nobreak
\medskip

 Above, Murre describes \G's feeling about the theory of the Picard
variety: it was hampered by its developers' narrow vision.  This section
explains how \G's broader vision led to clarifying and settling a number
of issues.  Primarily, we focus on the two major issues: Behavior in a
Family and Completeness of the Characteristic System.  In addition, we
consider some other issues mentioned earlier, especially Poincar\'e
divisors and the Albanese variety.   And we consider some other ways that
other mathematicians enhanced \G's theory between 1962 and 1973,
especially ways of generalizing his main existence theorem.
 For more discussion of those issues and some discussion of a lot of
other issues of the same sort, please see \cite{FGA}, \cite{FAG},
\cite{BLRnm}, and \cite{SGA6}.

As noted in the Introduction, \G\ explained the behavior of the Picard
schemes of the members of a family as compatibility with base
change.  More precisely, if the functor $\Pic_{(X/S)\fpqc}$ is
representable by an $S$-scheme $\pic_{X/S}$, then for any $S$-scheme
$S'$, the functor $\Pic_{(X\x S'/S')\fpqc}$ is representable by the
$S'$-scheme $\pic_{X/S}\x_SS'$.  In particular, if $S'$ is the spectrum
of the residue field $k_s$ of  $s\in S$, then the Picard scheme
of the fiber  $X_s$ of $X/S$ is just the fiber of $\pic_{X/S}$.

Compatibility holds for this reason.  For any $S'$-scheme $T$, the
equation
$$\Pic_{X\x S'/S'}(T) = \Pic_{X/S}(T)$$
 results from the definitions, because $(X\x_SS')\x_{S'}T=X\x_ST$.  So the
equation
$$\Pic_{(X\x S'/S')\fpqc}(T) = \Pic_{(X/S)\fpqc}(T)$$
 follows, because a map of $S'$-schemes $T'\to T$ is a covering if and
only if it is a covering when viewed as a map of $S$-schemes.  However,
the equation
$$\bigl(\pic_{X/S}\x_SS'\bigr)(T) = \pic_{X/S}(T)$$
 holds, because the structure map $T\to S'$ is fixed.  Since the right
sides of the last two equations are equal, so are their left sides, as
desired.

 {\bf Until otherwise said} near the end of the section, assume that $S$ is the
spectrum of an algebraically closed field $k$ and that $X$ is an
integral and projective $S$-scheme.  As is common, write ``$k$-scheme,''
$\IDiv_{X/k}$, etc.\ for ``$S$-scheme,'' $\IDiv_{X/S}$, etc.

In order to complete the discussion in Section 2 of the algebraic
proofs of the Theorem of Completeness of the Characteristic System and of
the Fundamental Theorem of Irregular Surfaces, we must discuss what's
called%
 \footnote{In 1947, Zariski \cite{Zsp} introduced and studied $\go m/\go m^2$
for a variety in any characteristic, and called it the ``local vector
space.''  In Weil's math review of Zarisk's paper, Weil wrote: ``the dual
vector-space \ldots\,seems to deserve to be called the
`tangent vector-space'.''}
  the {\it Zariski tangent space} $\bb T_z(Z)$ to a $k$-scheme $Z$ at a
{\it rational point} $z$, a point whose residue field $k_z$ is $k$.  Let
$\go m$ be the maximal ideal, and set $\bb T_z(Z):=\Hom_k(\go m/\go
m^2,\,k)$.

Then $\bb T_z(Z)$ can be viewed as the vector space of $k$-derivations
$\delta\:\mc O_{Z,z}\to k$; indeed, $\delta(\go m^2)=0$, and so $\delta$
corresponds to a linear map $\go m/\go m^2\to k$.  Let $k_\ve$ be the
ring of {\it dual numbers,} the ring obtained by adjoining an element
$\ve$ with $\ve^2=0$.  Then $\delta$ corresponds to the map of
$k$-algebras $u\:\mc O_{Z,z}\to k_\ve$ given by $u(a):=\?a+\delta(a)\ve$
where $\?a\in k$ is the residue of $a$.  Finally, let $S_\ve$ be the
spectrum of $k_\ve$; it is the {\it free tangent vector}.  Then $u$
corresponds to a $k$-map $S_\ve\to Z$, whose image is supported at $z$.
Denote the set of such $k$-maps by $h_Z(S_\ve)_z$.  Then in summary
$\bb T_z(Z)=h_Z(S_\ve)_z$.

Often, if  $Z$ represents a given functor $H$, that is $h_Z=H$,
 then we can work out a useful description of $h_Z(S_\ve)_z$ by viewing
it as the subset of $H(S_\ve)$ of elements whose image in $H(S)$ is $z$.
For example, say $Z=\hilb_{X/k}$ and $z\in Z$ represents $Y\subset X$.
Then $h_Z(S_\ve)_z$ is the set of $S_\ve$-flat closed subschemes of
$X\x_kS_\ve$ whose fiber over $S$ is $Y$.  Say the ideal of $Y$ is $\mc
I_Y$. Working it out, \G\ \cite[\pp.221-21--23]{FGA}
found $h_Z(S_\ve)_z=H^0(\mc N_Y)$ where $\mc N_Y:=\sHom(\mc I_Y/\mc
I_Y^2,\ \mc O_Y)$.  If $Y$ is a
Cartier divisor $D$, then $\mc N_D$ is invertible on $D$, and
\begin{equation}\label{eqchls}
 \bb T_z(\IDiv_{X/k})=H^0(\mc N_D).
\end{equation}

Let $\Lambda$ parameterize a system of divisors on $X$ including $D$,
and say $\lambda\in \Lambda$ represents $D$; in other words, there is a
map $\Lambda\to \IDiv_{X/k}$ carrying $\lambda$ to $z$.  It induces a
map of vector spaces $\theta\:\bb T_\lambda(\Lambda)\to \bb
T_z(\IDiv_{X/k})$.  If $D$ is integral, then, owing to \eqref{eqchls},
the image of $\theta$ defines a linear system on $D$, the storied {\it
characteristic linear system}.  When is it complete?  More generally,
for any $D$, when is $\theta$ surjective?

Each version of the Theorem of Completeness provides conditions
guaranteeing the existence of a $\Lambda$ that is smooth at $\lambda$
and for which $\theta$ is surjective.  But then $\IDiv_{X/k}$ is smooth
at $z$, owing to a simple general observation \cite[\p.305]{CCS}.  Thus
the conditions in question just guarantee that $\Lambda := \IDiv_{X/k}$
is smooth at $z$.

Some conditions are necessary.  Indeed, in 1943 Severi's student, Guido
Zappa, found a smooth complex surface $X$ such that $\IDiv_{X/k}$ has an
isolated point $z$ with $\dim\bb T_z(\IDiv_{X/k})=1$; so $\IDiv_{X/k}$
has nilpotents; for details, please see \cite[\pp.155--156]{Mcs} or
\cite[\p.285]{FAG}.  Commenting, \G\ \cite[\pp.221-24]{FGA} wrote that
this example ``shows in a particularly striking way how varieties with
nilpotents are needed to understand the phenomena of the most classical
theory of surfaces.''

\G\ then gave an enlightening proof of Kodaira's 1956
 version \cite{Kcls}
of Completeness.  As to Kodaira's own proof, Kodaira and Spencer
\cite[\p.477]{KSccs} said that it's ``based essentially on the theory of
harmonic differential forms'' [so not algebraic]; it's ``indirect and
does not reveal the real nature of the theorem.''

Kodaira proved that, {\it if $X$ and $D$ are smooth and complex, and\/%
 \footnote{It is now common to let $\mc O_X(D)$ stand for $\mc I_D^{-1}$,
but that practice is not followed here.}
  if $h^1(\mc I_D^{-1})=0$, then $\IDiv_{X/k}$ is smooth at $z$, where
$\mc I_D$ is the ideal of $D$.}  For
example, $h^1(\mc I_D^{-1})=0$ by Serre's computation if $D$ is a
hypersurface section of large degree.  In 1904, Enriques studied the
case that $X$ is a surface and $D$ is {\it regular\/} and {\it
nonspecial}, meaning%
 \footnote{\label{reg}At first, ``regular'' alone was used to mean
$h^1(\mc I_D^{-1})=0$ and $h^2(\mc I_D^{-1})=0$.}
  $h^1(\mc I_D^{-1})=h^2(\mc I_D^{-1})$ and $h^2(\mc
I_D^{-1})=0$.  Thus then Completeness holds.

\G\ proved Kodaira's theorem for any $X$ and $D$ as follows.  Let $\mc
I$ denote the ideal of the universal divisor on $X\x_k\IDiv_{X/k}$.
Then $\mc I^{-1}$ is invertible.  So it defines a map
$\balpha_{X/k}\:\IDiv_{X/k}\to \pic_{X/k}$, called the {\it Abel map}.
Assume $H^1(\mc I_D^{-1})=0$.  Then $\balpha_{X/k}$ is smooth at $z$;
see below.  Hence $\IDiv_{X/k}$ is smooth at $z$ if and only if
$\pic_{X/k}$ is smooth at $\balpha_{X/k}(z)$, or equivalently by
translation, everywhere.  In characteristic zero, $\pic_{X/k}$ is smooth
by Cartier's Theorem \cite[\p.167]{Mcs}.  Thus $\IDiv_{X/k}$ is smooth
at $z$ in characteristic zero.

In positive characteristic, $\smash{\pic_{X/k}}$ can be nonreduced
everywhere even if $X$ is a smooth surface; see below.  If so and
$h^1(\mc I_D^{-1})=0$, then $\IDiv_{X/k}$ is nonreduced at $z$. Thus
Completeness fails, even if $D$ is a general hypersurface section of
large degree.

Since $k$ is algebraically closed, $X$ has a rational point, so that
$X\to S$ has a section.  Set $P:=\pic_{X/k}$.  Then $\Pic_{X/k}$ is
representable by $P$.  So $X\x_kP$ carries an invertible sheaf $\mc P$,
called a {\it Poincar\'e sheaf}, whose class modulo $\Pic(P)$ is
universal.  In particular, there is an invertible sheaf $\mc M$ on $P$
such that $\mc I^{-1}$ is isomorphic to the tensor product on
$X\x_k\IDiv_{X/k}$ of the pullbacks of $\mc P$ and $\mc M$.  Then
$\LS_{\mc P/X\x P/P}$ is representable, on the one hand, by
$\IDiv_{X/k}$ regarded as a $P$-scheme via $\smash{\balpha_{X/k}}$, and
on the other, by $\I P(\mc Q)$ for some coherent sheaf $\mc Q$ on $P$.
So $\IDiv_{X/k}$ and $\I P(\mc Q)$ are canonically isomorphic
$P$-schemes.  If $H^1(\mc I_D^{-1})=0$, then $\mc Q$ is free at
$\balpha_{X/k}(z)$, and so $\balpha_{X/k}$ is smooth at $z$, as desired.

\G\ \cite[\pp.236-16]{FGA} asserted $\bb T_0(\pic_{X/k})=H^1(\mc
O_X)$; for proofs, please see \cite[\pp.163--164]{Mcs} and
\cite[\pp.281--282]{FAG}.  So $\dim \pic_{X/k}\le h^1(\mc O_X)$, with
equality if and only if $\pic_{X/k}$ is smooth.  That result is part of
\G's contribution to the proof of the Fundamental Theorem of Irregular
Surfaces.  In the examples of Igusa and Serre recalled in Section 2, we
have $\dim \pic_{X/k} < h^1(\mc O_X)$; hence, $\pic_{X/k}$ is not smooth
at 0, so nonreduced everywhere.

\G\ noted smoothness holds if $H^2(\mc O_X)=0$, owing to the
Infinitesimal Criterion for Smoothness and a well-known computation; for
details, please see \cite[\pp.285--286]{FAG}.  For example, if $X$ is a
curve, then $\pic_{X/k}$ is smooth, so of dimension $g$ where
$g:=h^1(O_X)$.  In positive characteristic, Mumford
\cite[\pp.193--198]{Mcs} proved  $\pic_{X/k}$ is smooth if and only
if  Serre's Bockstein operations \cite[\p.505]{ScpI} all vanish.

\G\ did not consider the other versions of Completeness, but his work
does provide a basis for proving them algebraically.  First consider
Severi's 1921 version\label{S21C} \eqref{eqccs} on \p.\pageref{eqccs}.
Given an invertible sheaf $\mc L$ on $X$, set $e(\mc L):=\chi(\mc
L)-1-h^2(\mc L)$.  Call $\mc L$ {\it arithmetically effective} if $e(\mc
L)\ge0$\label{ae}.  Vary $\mc L$.  Then $\chi(\mc L)$ is locally
constant, and $h^2(\mc L)$ is upper semi-continuous.  So $e(\mc L)$ is
lower semi-continuous.

Hence there is an open subscheme of $\pic_{X/k}$, say $\I P_{\rm ae}$,
that parameterizes the arithmetically effective $\mc L$ on $X$.  Set $\I
D_{\rm ae}:=\balpha_{X/k}^{-1}\I P_{\rm ae}$.  Assume $\dim X=2$. Then
$\I D_{\rm ae}$ surjects onto $\I P_{\rm ae}$, since over a point
representing an $\mc L$, the fiber has dimension $e(\mc L)+h^1(\mc L)$,
which is nonnegative.  In these terms, a refined Version \eqref{eqccs}
says that, {\it if $\pic_{X/k}$ is smooth too, then $\I D_{\rm ae}$ is smooth
on a dense open subset $U$.}

Since $h^0(\mc L)$ is upper semi-continuous in $\mc L$, there is an open
subscheme $V\subset\I P_{\rm ae}$ that parameterizes the $\mc L$ where
$h^0(\mc L)$ has a local minimum.  Set $U:=\balpha_{X/k}^{-1}V$.  Recall
that $\pic_{X/k}$ carries a coherent sheaf $\mc Q$ such that $\I P(\mc
Q)=\IDiv_{X/k}$.  Then $V$ is precisely the set of points of $\I P_{\rm
ae}$ at which the rank of $\mc Q$ has a local minimum.  Suppose
$\pic_{X/k}$ is smooth.  Then the restriction $\mc Q|V$ is locally free.
Hence $U\to V$ is smooth.  So $U$ is smooth.  Thus we have refined and proved
Severi's version \eqref{eqccs}.

In 1944, Severi discovered another condition on $D$ for Completeness to
hold if $\pic_{X/k}$ is smooth.  The condition requires $D$ to be {\it
semi-regular\/}; namely, in the standard long exact sequence of
cohomology
 \newbox\partbox\setbox\partbox=\hbox{$\partial^0$}
 \def\pto#1{\xto{\hbox to \wd\partbox{\scriptsize\hss$#1$\hss}}}
\begin{equation}\label{alles}
\begin{split}
 0&\To H^0(\mc O_X) \To H^0(\mc I_D^{-1}) \To H^0(\mc N_D) \\
  &\pto{\partial^0}  H^1(\mc O_X) \To  H^1(\mc I_D^{-1})
  \pto{u} H^1(\mc N_D)\pto{\partial^1}  H^2(\mc O_X),
\end{split}
\end{equation}
the map $u$ is 0, or equivalently $\partial^1$ is injective.
In particular, $D$ is semi-regular if either $ H^1(\mc I_D^{-1})=0$ or
$H^1(\mc N_D)=0$.

Severi worked with an integral $D$ on a smooth surface $X$, and he
formulated the condition in its dual form: the restriction
$H^0(\Omega^2_X)\to H^0(\Omega^2_X|D)$ is surjective; in other words,
the canonical system on $X$ cuts out a complete system on $D$.  In 1959,
Kodaira and Spencer \cite[\p.481]{KSccs} reformulated Severi's condition
as $u=0$ in any dimension.  Then they proved that, {\it in the complex
analytic case, if $X$ and $D$ are smooth and if $u=0$, then
$\IDiv_{X/k}$ is smooth at $z$.}

 \G\ did not consider semi-regularity per se, but he
\cite[\pp.221-23]{FGA} did observe that $H^1(\mc N_D)$ houses the
obstruction to deforming $D$ in $X$.  Thus%
 \footnote{Perhaps not surprisingly, the condition $H^1(\mc N_D)=0$ is
related to the flaw in the constructions of a good algebraic system made
in 1904 by Enriques and in 1905 by Severi.  In 1934, Zariski
\cite[\p.100]{Zas} noted that both constructions rely on a certain
assumption and that Severi's 1921 ``criticism is to the effect that the
available algebro-geometric proof of this assumption fails if'' $H^1(\mc
N_D)\neq0$.}
 {\it if $H^1(\mc N_D)=0$, then $\IDiv_{X/k}$ is smooth at $z$ in any
characteristic} whether $X$ and $D$ are smooth or not.  For example, if
$X$ is a curve, then $\IDiv_{X/k}$ is smooth everywhere; however,
$\balpha_{X/k}$ is not smooth at $z$ if $\deg D<g$ where $g:=h^1(\mc
O_X)$, since $\dim \IDiv_{X/k}=\deg D$ by \eqref{eqchls} and $\dim
\pic_{X/k}=g$ as noted above.

Mumford \cite[\pp.157--159]{Mcs} explicitly computed the obstruction to
deforming $D$, as well as its image under $\partial^1$.  He proved that
this image vanishes in characteristic 0 using an exponential. Therefore,
if $\partial^1$ is injective, then $\IDiv_{X/k}$ is smooth at $z$.
Cartier's Theorem is not involved, but recovered.  Thus in 1966 Mumford
gave the first algebraic proof that semi-regularity yields Completeness
in characteristic 0.

In 1973, I  \cite{CCS} gave another algebraic proof, yielding
a more refined statement: {\it assume $\pic_{X/k}$ is smooth;
then $\IDiv_{X/k}$ is smooth at $z$ of dimension $\rho$ where
\begin{equation*}\label{eqR}
\rho:= h^1(\mc O_X)-1 + h^0(\mc I_D^{-1}) - h^1(\mc I_D^{-1})
\end{equation*}
if and only if $D$ is semi-regular}.  My proof\,%
 \footnote{The proof works over any Noetherian $S$, and yields this more
general result: {\it let $D$ be a divisor on an (integral) geometric
fiber of $X/S$, and assume $\smash{\pic_{X/S}}$ is smooth; then
$\smash{\IDiv_{X/S}}$ is smooth of relative dimension $\rho$ at the
point representing $D$ if and only if $D$ is semi-regular}.}
 does not use obstruction theory, but a short formal analysis,
essentially due to George Kempf, of the scheme $\I P(\mc Q)$ above.

In passing, set $R:=\dim_z\IDiv_{X/k}$ and note that \eqref{eqchls}
yields $R\le h^0(\mc N_D)$, with equality if and only if $\IDiv_{X/k}$
is smooth at $z$.  Also, \eqref{alles} yields $\rho \le h^0(\mc N_D)$,
with equality if and only if $D$ is semi-regular.  Thus if $\IDiv_{X/k}$
is smooth at $z$, then $D$ is semi-regular if and only if $R=\rho$.

Generalizing more of Section 2, set $\delta:= \dim\Ck(\partial^0)$ and
$q:=h^1(\mc O_X)$.  Then \eqref{alles} yields $\delta\le q$, with
equality if $H^1(\mc I_D^{-1})=0$.  As $H^1(\mc I_D^{-1})=0$ if $D$ is a
hypersurface section of large degree, we have generalized Castelnuovo's
result (1).  Next, set $r:=h^0(\mc I_D^{-1})-1$.  Then \eqref{alles}
yields $h^0(\mc N_D)=r+\delta$.  Hence $R\le r+\delta$, with equality if
and only if $\IDiv_{X/k}$ is smooth at $z$.  Thus we have generalized
Severi's result \eqref{eqCpltnss}.  Finally, {\it if $X$ is a surface,
$\pic_{X/k}$ is smooth and $z$ lies in the open subset $U\subset \I
D_{\rm ae}$, then $\IDiv_{X/k}$ is smooth at $z$, and so
$R= r+\delta$}, just as Severi discovered.

 \G\ \cite[\pp.2-12]{FGA} proved the following basic properties of the
connected component of $0$ in $\pic_{X/k}$, denoted $\pic_{X/k}^0$.  It
is open and closed.  It is irreducible.  Forming it commutes with base
change.  It is {\it quasi-projective\/}; that is, it is an open
subscheme of a projective $k$-scheme.  Moreover, it is projective if $X$
is normal.  Of course, the {\it Picard variety\/} of $X$ is the set of
points of $\smash{\pic_{X/k}^0}$ with coordinates in a given universal
domain.  If $X$ is a curve, then $\pic_{X/k}^0$ is its {\it generalized
Jacobian}.

Define a {\it Poincar\'e divisor\/} to be a Cartier divisor $\Delta$ on
$X\x P$, where $P$ is some connected component of $\pic_{X/k}$, such
that $\Delta$ yields a section $P\to\IDiv_{X/k}$.  Such $\Delta$ abound,
as is shown next, generalizing Mattuck's work \label{Mat} mentioned on
\p.\pageref{Mat1}.

Given any connected component $P'$ of $\pic_{X/k}$, notice it's a
translate of $\pic_{X/k}^0$, so quasi-projective.  Hence there's an $n$
such that, for any invertible sheaf $\mc L$ on $X$ represented by a
point of $P'$, we have $h^0(\mc L(n))>\dim \pic_{X/k}$ and $h^1(\mc
L(n))=0$ where $\mc L(n)$ is the $n$th twist by
$\mc O_X(1)$.  Let $P$ be the translate of $P'$ defined by $\mc O_X(n)$.
Since $P$ is  quasi-projective, it has a universal sheaf $\mc O_P(1)$. 

Recall $\pic_{X/k}$ carries a coherent sheaf $\mc Q$ such that $\I P(\mc
Q)=\IDiv_{X/k}$.  Notice the restriction $\mc Q|P$ is locally free of
rank $h^0(\mc L(n))$, so of rank more than $\dim P$.  Take $m$
so that $\sHom(\mc Q|P,\,\mc O_P)(m)$ is generated by its global
sections, so by finitely many.  Then a general linear combination of
the latter vanishes nowhere on $P$ by a well-known lemma, \cite[\p.426]{Avbec}
or \cite[\p.148]{Mcs}, due to Serre.  So there's a surjection $\mc
Q|P\onto \mc O_P(m)$.  It defines a section $P\to \I P(\mc Q)$, and so a
Poincar\'e divisor $\Delta$.

Suppose also that $X$ is a curve.  Set $g:=h^1(\mc O_X)$ and
recall $g=\dim \pic_{X/k}$.  Take $P$ to be any connected component of
$\pic_{X/k}$ that parameterizes invertible sheaves $\mc L$ on $X$ of
$\deg\mc L>2g-1$.  Then $h^0(\mc L)> g$ and $h^1(\mc L)=0$.  So similarly
there is a section $P\to \I P(\mc Q)$, and so a Poincar\'e divisor
$\Delta$.

Here's an introduction\label{GAlb}%
 \footnote{Doubtless, \G\ had something similar in mind when he
\cite[\p.232-14]{FGA} wrote, ``the theory of Abelian varieties, and more
generally of Abelian schemes, becomes much simpler once we have a
general theory of the Picard scheme at our disposal.  In particular, the
theory of duality for Abelian schemes, and notably results like
Cartier's, thus become nearly formal (cf.\ for example the forthcoming
notes to the Mumford--Tate seminar at Harvard in the spring term of
1962).''  However, it seems that nothing like the present introduction
appears in Mumford's personal notes to the seminar, or has already
appeared in print.}
  to the scheme-theoretic theory of the Albanese variety.  Assume $X$ is
normal.  Then $\pic_{X/k}^0$ is projective.  Let $P$ denote its {\it
reduction}, namely, the subscheme defined by the nilradical of the
structure sheaf of $\smash{\pic_{X/k}^0}$.  It too is a {\it group
scheme}; that is, its $T$-points form a group for all $T$.  So $P$ is
smooth.  Call any such connected smooth projective group scheme an {\it
Abelian variety}.

If $X$ is an Abelian variety, then $\pic_{X/k}^0$ is already reduced.
Mumford \cite{Mgit} gave a proof on \pp.117--118, which he attributed to
\G\ on \p.115.  Then $\smash{\pic_{X/k}^0}$ is denoted by $\smash{\wh
X}$ or $X^*$, and called the {\it dual Abelian variety}.

In general, let $Y$ be another integral and projective $k$-scheme, and
fix rational points $x\in X$ and $y\in Y$.  Then a $k$-map $f\:Y\to P$
with $f(y)=0$ is defined by an invertible sheaf $\mc L$ on $X\x_kY$
whose restriction to $X\x_k y$ is $\mc O_X$.  At first, $\mc L$ is only
determined modulo $\Pic(Y)$, but normalize $\mc L$ as follows: restrict
it to $x\x_kY$, pull the restriction back to $X\x_kY$ via the
projection, and replace $\mc L$ by its tensor product with the inverse
of the pullback.

Let $Q$ be the reduction of $\pic_{Y/k}$.  By symmetry, $\mc L$
corresponds to a $k$-map $g\:X\to Q$ with $g(x)=0$.  Plainly, this
correspondence $f\longleftrightarrow g$ is {\it functorial\/}: given a
similar triple $(Z,\, z,\, R)$ and a $k$-map $h\:Z\to Y$ with $g(z)=y$,
the composition $fh\:Z\to P$ corresponds to the composition $h^*\circ
g\:X\to R$ where $h^*\:P\to R$ is the $k$-map induced by pullback of
invertible sheaves, which is a group homomorphism.

For a moment, take $Y:=P$ and $y:=0$ and $f:=1_p$.  Since $P$ is an
Abelian variety, $Q$ is its dual $P^*$.  Set $A:=P^*=Q$ and $a:=g$.
Then $A$ is called the {\it Albanese variety\/} of $X$, and there is a
canonical $k$-map $a\:X\to A$.

The map $a\:X\to A$ is the {\it universal example\/} of a map $g\:X\to
Q$ where $Q$ is the reduction of $\pic_{Y/k}$ for some integral and
projective $k$-scheme $Y$; that is, any such $g$ factors uniquely
through $a$.  Here's why.  Say $g$ corresponds to $h\:Y\to P$.  Then by
functoriality, $1_P\circ h$ corresponds to $h^*\circ a$.

By definition, $A^*$ is the Albanese of $P$.  Moreover, the canonical
map $p\:P\to A^*$ is an isomorphism, since by ``abstract nonsense,'' a
universal example is determined up to unique isomorphism, and $1_P\:P\to
P$ is trivially another universal example.  In fact, $p^{-1}$ is just
the map $a^*$ induced by the canonical map $a\:X\to A$, because by
functoriality, $1_A\circ a$ corresponds to $a^*\circ p$.
Thus $A$ and $P$ are dual to each other.

Suppose $X$ is an Abelian variety;  take $x:=0$.  Mumford \cite[\p.125]{Mav}
constructed $X^*$ as a quotient of $X$ by a finite subgroup; so $X$ and
$X^*$ are isogenous.  Mumford \cite[Cor., \p.43,\,132]{Mav} proved that
$a:X\to X^{**}$ is an isomorphism of groups and of schemes.  Hence, for
any $X$, the map $a\:X\to A$ is the {\it universal example\/} of a map
$X\to B$ where $B$ is an Abelian variety, because $B=Q$ if $Y=B^*$.

Suppose finally that $X$ is a smooth curve of genus $g>0$.  Then $X$ is
a component of $\IDiv_{X/k}$.  So the Abel map restricts to a $k$-map
$X\to \pic_{X/k}$.  Its image lies in the connected component
parameterizing the sheaves $\mc L$ of degree 1.  Fix an $\mc L$.
Translating by $\mc L^{-1}$ yields a map $\alpha\:X\to P$.  It is an
embedding by general principles, since its fibers are finite and $X=\I
P(\mc Q)$ for some coherent sheaf $\mc Q$ on $P$.  It is proved (in a
more general form) in \cite[Thm.\,2.1, \p.595]{EGK} that
$\alpha^*\:P^*\to P$ is an isomorphism, which is independent of the
choice of $\mc L$.

 {\bf To end} this article, let's consider some important ways in which
\G\ and others generalized the existence theorem culminating Section~4.
First, let $k$ be an arbitrary field.  On \pp.232-15--17 in \cite{FGA},
\G\ outlined a construction of $\pic_{X/k}$ for any projective
$k$-scheme $X$.  He used that earlier theorem plus a method of {\it
relative representability}, by which $\pic_{X/k}$ is constructed from
$\pic_{X'/k}$ for a suitable surjective $k$-map $X'\onto X$.  The method
employs two main tools: nonflat descent and Oort d\'evissage.  The
former refers to descent along maps not required to be flat; however,
key objects are flat.  The second tool was introduced by Oort in
\cite{OsP} to construct $\pic_{X/k}$ from $\pic_{X'/k}$ where $X'$ is
the reduction of $X$.

On \p.232-17 in \cite{FGA}, \G, in effect, made two conjectures:
first, $\pic_{X/k}$ exists for any proper $k$-scheme $X$;
second, given any surjective $k$-map between proper $k$-schemes, the
induced map on Picard schemes is affine.

The first conjecture was proved in 1964 by Murre \cite{Mcf}, who thanked
\G\ for help.  However, instead of using relative representability,
Murre identified seven conditions that are necessary and sufficient for
the representability of a functor from schemes over a field to Abelian
groups.  Then he checked the seven for $\Pic_{(X/k)\fpqc}$.

 {\bf From now on}, assume $S$ is Noetherian and $X$ is a flat and proper
$S$-scheme.

Murre \cite[\p.5]{Mcf} said that \G\ too proved the first conjecture.  In
1965, Murre \cite{Mru} sketched \G's proof of the following key
intermediate result: {\it let $\mc F$ be a coherent sheaf on $X$, and
$S_{\mc F}$ the functor of all $S$-schemes $T$ such that the pullback
$\mc F_T$ is $T$-flat; then $S_{\mc F}$ is representable by an
unramified $S$-scheme of finite type.}  The proof involves identifying
and checking eight conditions that are necessary and sufficient for
representability by a scheme of the desired sort.

In 1966, Raynaud \cite[Exp.~XII]{SGA6} gave \G's proof of another
key intermediate result: {\it assume $S$ is integral and let $X'\onto X$
be a surjective map of proper $S$-schemes; then there's a nonempty open
subscheme $V\subset S$ such that $\pic_{X'\x V/V}$ and $\pic_{X\x V/V}$
exist, and the induced map between them is quasi-affine.}  The proof
does indeed involve suitably general versions of nonflat descent and
Oort d\'evissage.  As corollaries, that result yields \G's two
conjectures.

If the geometric fibers of $X/S$ are not all integral, then $\pic_{X/S}$
need not exist.  On \p.236-01 in \cite{FGA}, \G\ described an example of
Mumford's; one geometric fiber is integral, but another is a pair of
conjugate lines.  On the other hand, on \p.viii in \cite{Mcs}, Mumford 
asserted this theorem: {\it Assume $X/S$ is projective, and its
geometric fibers are reduced and connected; assume the irreducible
components of its ordinary fibers are geometrically irreducible; then
$\pic_{X/S}$ exists.}  He said the proof is like the one on
\pp.133--149, involving his theory of independent 0-cycles.

On \p.236-01, \G\ attributed a slightly different theorem to Mumford,
and referred to the Mumford--Tate seminar.  Mumford's seminar notes
contain a precise statement of the theorem and a rough sketch of the
proof.  However, he crossed out the hypothesis that the geometric fibers
are connected, and made the weaker assumption that the ordinary fibers
are connected.

On \p.236-13, \G\ wrote that ``it is not ruled out that $\pic_{X/S}$
exists'' whenever%
 \footnote{This condition holds if the geometric fibers of $X/S$ are
integral by \cite[Prp.\ (7.8.6), \p.74]{EGAIII}.  For more about its
significance when $S$ is the spectrum of a discrete valuation
ring, please see \cite{Rsfp}.}
 the direct image of $\mc O_{X\x T}$ is $\mc O_{T}$ for any $T$.  ``At
least, this statement is proved for analytic spaces when $X/S$ is also
projective.''  Mumford's example shows the statement is false for
schemes.  Michael Artin's work shows it holds for {\it algebraic
spaces}, which he introduced in 1968 in \cite{Aif}.  They are formed by
gluing together schemes along open subsets that are isomorphic locally
in the \'etale topology.  Over $\bb C$, these open sets are locally
analytically isomorphic; so a separated algebraic space is a kind of
complex analytic space.

In 1969, Artin \cite{Afm}, inspired by Grothendieck and Murre, found
five conditions on a functor that are necessary and sufficient for it to
be representable by a well-behaved algebraic space.  A key new
ingredient is Artin's Approximation Theorem; it facilitates the
passage from formal power series to polynomials.  By checking that the
conditions hold if the direct image of $\mc O_{X\x T}$ is always $\mc
O_{T}$, Artin \cite[Thm.\,7.3, \p.67]{Afm} proved%
 \footnote{In fact, he proved a more general theorem, in which $S$ and
$X$ are algebraic spaces.  That theorem and Grothendieck's theorem in
Section~4 are the two main representability theorems for the Picard
functor.  Grothendieck used projective methods.  Artin's work has a very
different flavor.  Moreover, it yields a major improvement of Murre's
representability theorem stated above, and it yields the
representability of the Hilbert functor and related functors in
algebraic spaces.}
  $\pic_{X/S}$ exists as an algebraic space, a magnificent achievement.
Also, he \cite[Lem.\,4.2, \p.43]{Afm} proved that, if $S$ is the
spectrum of a field, then $\pic_{X/S}$ is a scheme.  Thus he obtained a
third proof of \G's first conjecture.

As $S$ and $X$ are schemes, so are the fibers of $X/S$.  Hence their
Picard schemes exist.  Furthermore, if the direct image of $\mc O_{X\x
T}$ is always $\mc O_{T}$, then these Picard schemes form a family; its
total space $\pic_{X/S}$ is an algebraic space, but need not be a
scheme.  Thus Artin proved the definitive statement explaining the
behavior of the Picard schemes of the members of a family.%
 \footnote{As the Picard varieties in the family are the points of the
component $\pic_{X_s/k_s}^0$ for $s\in S$, their behavior is explained
by an open subspace $\pic_{X/S}^0$ of $\pic_{X/S}$ whose fibers are the
$\pic_{X_s/k_s}^0$.  Such a  $\pic_{X/S}^0$ is observed
in \cite[\p.233]{BLRnm} to exist when  $\pic_{X/S}$ is  $S$-smooth along
the 0-section.}

\bigbreak

\end{document}